\documentclass[11pt]{amsart}   
\usepackage{Angelini-Knoll} 
\usepackage[margin=1.2in]{geometry} 

\DeclareMathOperator{\et}{\text{\'et}}

\DeclareMathOperator{\THH}{\rm THH}
\DeclareMathOperator{\HH}{\rm HH}
\DeclareMathOperator{\K}{\rm K}  
\DeclareMathOperator{\TC}{\rm TC}
\DeclareMathOperator{\TP}{\rm TP} 
\DeclareMathOperator{\Assoc}{\operatorname{Assoc}}

\usepackage{enumitem}

\usepackage[normalem]{ulem}	  
\makeatletter
\def\blfootnote{\gdef\@thefnmark{}\@footnotetext} 
\makeatother
\usepackage[colorlinks=true, linkcolor=blue, menucolor=blue]{hyperref}
\usepackage{amsrefs}
\numberwithin{equation}{section}

\title{Detecting $\beta$ elements in iterated algebraic K-theory}
\author{Gabriel Angelini-Knoll} 
\address{Department of Mathematics, Michigan State University, East Lansing, Michigan 48824
}
\address{\emph{Current address:} Department of Mathematics, Institut Galil\'ee,
Universit\'e Sorbonne Paris Nord, 99 av. JB Clément, FR-93430 Villetaneuse, France
}

\begin{document} 
\blfootnote{\textup{2010} \textit{Mathematics Subject Classification}:
Primary 19D55; Secondary 55P42}
\begin{abstract} 
The Lichtenbaum--Quillen conjecture (LQC) relates special values of zeta functions to algebraic K-theory groups. The Ausoni--Rognes red-shift conjectures generalize the LQC to higher chromatic heights in a precise sense. In this paper, we propose an alternate generalization of the LQC to higher chromatic heights and give evidence for it at height two. In particular, if the $n$-th Greek letter family is detected by a commutative ring spectrum $R$, then we conjecture that the $n+1$-st Greek letter family will be detected by the algebraic K-theory of $R$. We prove this in the case $n=1$ for $R=\K(\mathbb{F}_q)$ modulo $(p,v_1)$ where $p\ge 5$ and $q=\ell^k$ is a prime power generator of the units in $\mathbb{Z}/p^2\mathbb{Z}$. In particular, we prove that  the commutative ring spectrum $\K(\K(\mathbb{F}_q))$ detects the part of the 
$p$-primary $\beta$-family that survives mod $(p,v_1)$. The method of proof also implies that these $\beta$ elements are detected in iterated algebraic K-theory of the integers. Consequently, one may relate iterated algebraic K-theory groups of the integers to integral modular forms satisfying certain congruences.   
\end{abstract}
\maketitle

\tableofcontents
\section{Introduction}
Following Waldhausen \cite[\S 4]{MR764579}, the famous Lichtenbaum--Quillen conjecture \cite{MR0422392,MR0406981} states that the map
\begin{align}\label{LQmap} \K_n(A;\bZ/p\bZ)\to \K_n^{\et}(A; \bZ/ p\bZ) \end{align}
from algebraic K-theory to \'etale K-theory is an isomorphism for $n$ sufficiently large for suitable rings $A$ when $p>2$ is a prime. Here by suitable ring $A$, we mean that $A$ is a regular Noetherian commutative ring with $A=A[p^{-1}]$, which is finitely generated as a $\mathbb{Z}$-algebra \cite[\S 9 p.175]{MR0406981}.
Since algebraic K-theory satisfies Nisnevich descent and \'etale K-theory satisfies \'etale descent, the question can be translated into the question of whether the map from motivic cohomology to \'etale cohomology is an isomorphism in a range. In this way, the conjecture was resolved by work of M. Rost and V. Voevodsky \cite{Voe11} as a consequence of their proof of the Bloch--Kato conjecture (see \cite{BKbook} for a complete proof and \cite[Historical Remark 4.4]{Kbook} for further discussion).

In \cite{MR826102}, R.W. Thomason showed  that $\K_n^{\et}(A; \bZ/ p \bZ)\cong \beta^{-1}\K_n(A; \bZ/ p\bZ)$ under the same conditions on $A$ where $\beta$ is the Bott element in $\K_{2}(A; \bZ/ p\bZ)$. From the perspective of homotopy theory, we may therefore view the map \eqref{LQmap} as the map on $\pi_n$ induced by the map of spectra 
\[ S/p \wedge \K(A) \rightarrow v_1^{-1}S/p \wedge \K(A) \] 
where $S/p$ is the cofiber of multiplication by $p$. 
Here the map $v_1\co \Sigma^{2p-2} S/p \to S/p$ is a $v_1$-self map, which means it is non-nilpotent, and we define $v_1^{-1}S/p$ as the homotopy colimit of the diagram
\[S/p \overset{v_1}{\longrightarrow} \Sigma^{-2p+2}S/p \overset{v_1}{\longrightarrow} \Sigma^{-4p+4}S/p \overset{v_1}{\longrightarrow} \ldots \]
of spectra. The effect of inverting the Bott element is the same as the effect of inverting $v_1$ by work of Snaith \cite{MR750689} as interpreted by Waldhausen \cite[\S 4]{MR764579}. 

The original motivation of the Lichtenbaum--Quillen conjecture was to relate algebraic K-theory groups to special values of zeta functions. For $A$ the ring of integers in a totally real number field $F$ and $p$ an odd prime, Wiles proved that quotients of \'etale cohomology groups of $A[1/p]$ recover special values of the Dedekind zeta function $\zeta_F$ \cite{MR1053488}. 
The Lichtenbaum--Quillen conjecture then gives a correspondence between algebraic K-theory groups and special values of Dedekind zeta functions, see \cite[Conjecture 2.4]{MR0422392}. Notably these special values correspond to the $v_1$-periodic part of $S/p _*\K(A)$ because they are detected in $v_1^{-1}S/p _*\K(A)$. 

Fix an odd prime $p$ throughout this section and fix a prime power $q$ that generates the units in $\mathbb{Z}/p^2\mathbb{Z}$.
As another specific example, consider the algebraic K-theory of finite fields $\mathbb{F}_q$. D. Quillen \cite{MR0315016} computed $\K_n(\mathbb{F}_q)$ for all $n$ and after localizing at $p$, there is an isomorphism
\[\K_{2(p-1)k-1}(\mathbb{F}_q;\bZ_{(p)})\cong \bZ/p^{\nu_p(k)+1}\bZ\] 
where $\nu_{p}(k)$ is the $p$-adic valuation of $k$. 
The order of the group $\K_{2s-1}(\mathbb{F}_q;\bZ_{(p)})$ corresponds exactly to the $p$-adic valuation of the denominator of $B_{s}/2s$ where $B_s$ is the $s$-th Bernoulli number. Recall that Bernoulli numbers are the coefficients in the Taylor series 
\[ \frac{x}{e^x-1}=\sum_{s\ge 0} B_s \frac{x^s}{s!} \]
and the special values of the Riemann zeta function satisfy $\zeta(-s)=(-1)^sB_s/(s+1)$ for $s\ge 0$. 

This example is intimately tied to stable homotopy theory as well. J. F. Adams showed that the image of the J-homomorphism from the homotopy groups of the stable orthogonal group to the stable homotopy groups of spheres is highly nontrivial and the classical Bott periodicity in the homotopy groups of the stable orthogonal group corresponds to periodicity in the stable homotopy groups of spheres \cite{MR0198470}. In fact, the $p$-local image of the J-homomorphism exactly corresponds to the image of the map $\pi_{2(p-1)k-1}S_{(p)}\to \K_{2(p-1)k-1}(\mathbb{F}_q;\bZ_{(p)})$ when $p$ is an odd prime. The image of J therefore bridges the fields of homotopy theory and number theory. The spectrum $H\mathbb{F}_q$ detects $v_0$-periodicity in the sense that $p^k=v_0^k$ is nontrivial in the image of the Hurewicz map $\pi_*S\to \pi_*H\mathbb{F}_q$. Therefore, we have observed an instance where algebraic K-theory of a spectrum that detects $v_0$-periodic elements detects $v_1$-periodic elements. One goal of this introduction is to formulate a precise conjecture about a higher chromatic height generalization of this phenomenon. The main theorem of this paper is evidence for this conjecture at a higher chromatic height. 

In chromatic stable homotopy theory, we study periodic families of elements in the homotopy groups of spheres. The first such family, due to J. F. Adams \cite{MR0198470} and H. Toda \cite{MR0111041}, is the $\alpha$-family, which consists of maps $\alpha_k$ defined as the composites 
\[ \xymatrix{ \alpha_k\co \Sigma^{(2p-2)k}S \ar[r]^{i_0} & \Sigma^{(2p-2)k}S/p \ar[r]^(.6){v_1^k} & S/p \ar[r]^{\delta_0}& \Sigma S }\]  
where $p$ is an odd prime. The elements $\alpha_k$ are $p$-torsion elements in the groups $\pi_{2(p-1)k-1}S$. These elements are in the image of the J-homomorphism at odd primes $p$ and as discussed earlier they are also detected in algebraic K-theory of finite fields of order $q$, under our hypotheses on $p$ and $q$. 
In particular, they correspond to certain special values of the Riemann zeta function. Now, consider the cofiber of the periodic self-map $v_1\co\Sigma^{2p-2}S/p\rightarrow S/p$ denoted $V(1)$. When $p\ge 5$, there exists a periodic self-map $v_2\co \Sigma^{2p^2-2}V(1)\to V(1)$ and there is an associated periodic family of elements in the homotopy groups of spheres called the $\beta$-family. In particular, L. Smith \cite{MR0275429} proved that the maps 
\[ \xymatrix{ \beta_k\co \Sigma^{(2p^2-2)k}S\ar[r]^{i_1i_0} & \Sigma^{(2p^2-2)k}V(1)\ar[r]^(.6){v_2^k} & V(1) \ar[r]^{\delta_0\delta_1} & \Sigma^{2p} S } \]
are nontrivial. 

In the language of chromatic homotopy theory, the $\alpha$-family is a periodic family of height one and the $\beta$-family is a periodic family of height two. There are a family of homology theories $K(n)_*$ called Morava K-theory which are useful for detecting periodicity of chromatic height $n$ in the homotopy groups of spheres. The coefficients of Morava K-theory are $K(n)_*\cong \mathbb{F}_p[v_n^{\pm 1}]$ for $n\ge 1$ and $K(0)_*$ is rational homology. We say a $p$-local finite cell $S$-module $V$ has type $n$ if the groups $K(n)_*V\neq 0$ and the groups $K(n-1)_*V$ vanish. By the celebrated periodicity theorem of
Hopkins-Smith \cite{MR1652975}, any $p$-local finite spectrum $V$ of type $n$  admits a periodic self map 
\[ v_n^m\co\Sigma^{(2p^n-2)m}V\rightarrow V. \]  
We can therefore define $v_n^{-1}V$ in the same way that we defined $v_1^{-1}S/p$. We can also construct the $n$-th Greek letter family by including into the bottom cell, iterating $v_n^m$ $k$-times, and then projecting onto the top cell. However, it is highly non-trivial to prove that Greek letter elements that are constructed in this way are actually nonzero.

The study of Greek letter family elements was significantly expanded by the groundbreaking work of Miller--Ravenel--Wilson \cite{MR0458423} using the chromatic spectral sequence
\begin{align}\label{chrom ss} E_1^{*,*}= \bigoplus_{i\ge 0} \Ext_{BP_*BP}^{*,*}(BP_*,v_i^{-1}BP_*/(p^{\infty},v_1^{\infty},\dots,v_{i-1}^{\infty}))\implies \Ext_{BP_*BP}^{*,*}(BP_*,BP_*) \end{align}
which converges to the input of the $BP$-Adams spectral sequence. 
If the class 
\begin{align*} 
v_n^k/p^{i_0}v_1^{i_1}\dots v_{n-1}^{i_{n-1}}\in \Ext_{BP_*BP}^0(BP_*,v_n^{-1}BP_*/(p^{\infty},v_1^{\infty},\dots,v_{n-1}^{\infty}))
\end{align*}
in the $E_1$-page of \eqref{chrom ss} is a permanent cycle in the chromatic spectral sequence, we will write 
\begin{align}\label{Greek letter}
\overline{\alpha}^{(n)}_{k/(i_{n-1},i_{n-2},\dots i_0)}\in \Ext_{BP_*BP}^{*,*}(BP_*,BP_*) 
\end{align}
for its image in the abutment of the chromatic spectral sequence. We will refer to the collection of all such elements for a fixed $n$ as the $n$-th divided (algebraic) Greek letter family  and when any of the $i_j$ for $0\le j\le n-1$ are $1$ we omit them from the notation. If the elements $\overline{\alpha}^{(n)}_k$ are permanent cycles in the $BP$-Adams spectral sequence, then we will write $\alpha^{(n)}_k$ for the corresponding elements in the stable homotopy groups of spheres and refer to the collection as the \emph{$n$-th Greek letter family}. The advantage of this approach is that the elements in the input of the chromatic spectral sequence always exist. The question of whether or not certain Greek letter elements exist in homotopy can then be approached by determining whether certain elements in the chromatic spectral sequence and the $BP$-Adams spectral sequence are permanent cycles. 

We will say that a ring spectrum $R$ \emph{detects} the $n$-th Greek letter family $\{\alpha^{(n)}_k\}$ in the homotopy groups of spheres if each element $\alpha^{(n)}_k$ is non-trivial in the image of the unit map 
\[ \pi_*S\lra \pi_*R.\] 
We conjecture the following higher chromatic height analogue of the Lichtenbaum--Quillen conjecture, which is in the same spirit as the red-shift conjectures of Ausoni--Rognes \cite{AR08}. For the following conjecture, suppose the $n$-th and the $n+1$-st Greek letter family are nontrivial elements in $\pi_*S$ for a given prime $p$. 
\begin{conj}\label{Greek red} 
If $R$ is a commutative ring spectrum that detects the $n$-th Greek letter family, then $\K(R)$ detects the $n+1$-st Greek letter family. 
\end{conj}
When $p$ is sufficiently small this conjecture is both more difficult to verify and more subtle, for example at the prime $p=3$ only the $\beta$ elements $\beta_i$ for $i\equiv 0,1,3,5,6\mod 9$ are know to exist by \cite{BP04}. 

We can now state the main theorem of this paper.  As discussed earlier, the spectrum $\K(\bF_q)_p$ detects the $\alpha$-family for $p\ge 3$. The main theorem of this paper is a proof of Conjecture \ref{Greek red} in the case $n=1$ where $R=\K(\bF_q)_p$, with the following caveat: Since the approach we take is to detect the $\beta$-family in the image of the unit map 
\[ S_*\to V(1)_*\K(\K(\mathbb{F}_q)_p)\]
we only detect the part of the $p$-primary $\beta$-family that is detected in $V(1)_*$.\footnote{Note that the elements $\beta_{pk}$ are not $p$-divisible by \cite[Theorem 2.12]{MR0458423} and they are $p$-torsion by construction, however they are $v_1$-divisible in $\pi_*S/p$ by \cite[Theorem A]{Oka75} and \cite[Theorem AII]{Oka76} and consequently they have trivial Hurewicz image in $V(1)_*$, see \cite{Shi14} for a survey.}

\begin{thm}
The commutative ring spectrum $\K(\K(\bF_q)_p)$ detects the $p$-primary $\beta$-family $\{\beta_k: k\not \equiv 0\mod p\}$ for all $p\ge 5$. 
\end{thm}
In particular, the method of proof also provides the following higher Lichtenbaum--Quillen-type result about iterated algebraic K-theory of the integers. 
\begin{cor}
The commutative ring spectrum $\K(\K(\mathbb{Z}))$ detects the $p$-primary $\beta$-family $\{\beta_k: k\not \equiv 0\mod p\}$ for all $p\ge 5$. 
\end{cor}
In \cite{MR2469520,MR2544384,MR1660325}, M. Behrens and G. Laures describe an explicit relationship between the $\beta$-family and certain integral modular forms satisfying certain congruences. From this point of view, our main result may be viewed as a higher chromatic height version of the Lichtenbaum--Quillen conjecture. It is therefore a step towards the larger program of understanding the arithmetic of commutative ring spectra. 

The $\beta$ elements $\beta_k$ that we detect only agree with the divided $\beta$-family elements $\beta_{k/i,j}$ of M. Behrens \cite{MR2469520} when $i=j=1$ and $k\not\equiv 0 \mod p$. The explicit modular forms associated to these beta elements are also computed in \cite{Lar19}. It would be desirable to detect the entire \emph{divided} $\beta$-family in iterated algebraic K-theory of finite fields to detect more information about integral modular forms. 

\subsection{Conventions}
Let $\Sp$ be the category of symmetric spectra in pointed simplicial sets with the positive flat stable model structure. Most of the results here can also be proven for other models of the stable homotopy category since they depend only on the homotopy category, but the proof relies on \cite{thhmay}, which uses this model structure. We write $\Comm \Sp$ for the category of commutative monoids in $\Sp$, which we call commutative ring spectra, and $\Assoc \Sp$ for the category of associative monoids in $\Sp$, which we call associative ring spectra. By a ring spectrum, we mean a homotopy associative monoid in $\Sp$.  

Throughout, we write $H_*(-)$ for $\pi_*(-\wedge H\bF_p)$; i.e, homology with $\bF_p$-coefficients and for a general spectrum $\mathsf{E}$ we write $\mathsf{E}_*X=\pi_*(X\wedge \mathsf{E})$. We also write $H_*(-;d)$ for homology of a chain complex with  differential $d$ and dually we write $H^*(-;d)$ for the cohomology of a cochain complex with differential $d$. We write $\otimes$ for tensor over $\mathbb{F}_p$ if not otherwise specified. We write $\HH_*(A)$ for the Hochschild homology over $\mathbb{F}_p$ of a bi-graded $\mathbb{F}_p$ algebra $A$. We also fix the convention that $E_R(x_1, \dots ,x_n)$ denotes the graded exterior algebra over a ring $R$ with indecomposable algebra generators $x_1,\dots ,x_n$, $P_R(y_1, \dots , y_n)$ denotes the graded polynomial algebra over a ring $R$ with indecomposable algebra generators $y_1,\dots ,y_n$, and $\Gamma_{R}(z_1,\dots z_n)$ denotes the graded divided power algebra with indecomposable algebra generators $z_1,\dots ,z_n$. When $R=\mathbb{F}_p$ we omit $R$ from the notation. 

A comodule $M$ over a Hopf algebroid $(\mathsf{E}_*,\mathsf{E}_*\mathsf{E})$  associated to a ring spectrum $\mathsf{E}$ will always be considered with left co-action 
\[ \psi_M^{\mathsf{E}} \colon \thinspace M\rightarrow \mathsf{E}_*\mathsf{E}\otimes_{\mathsf{E}_*} M \]
and we will simply write $\psi$ when the module $M$ and the Hopf algebroid $(\mathsf{E},\mathsf{E}_*\mathsf{E})$ is understood from the context. The main examples of interest are $\mathsf{E}=H\bF_p$, where $\mathsf{E}_*\mathsf{E}$ is the dual Steenrod algebra $\mathcal{A}_*$, and $\mathsf{E}=BP$. We write 
\[ \Delta_{\mathsf{E}} \colon \thinspace \mathsf{E}_*\mathsf{E} \rightarrow \mathsf{E}_*\mathsf{E} \otimes_{\mathsf{E}_*} \mathsf{E}_*\mathsf{E}\]
for the co-product of the Hopf-algebroid $\mathsf{E}_*\mathsf{E}$ or simply $\Delta$ when $\mathsf{E}$ is understood from the context. When $\mathsf{E}=H\bF_p$, this is the co-product in the dual Steenrod algebra $\mathcal{A}_*\cong P(\bar{\xi}_i \text{ }|\text{ }i\ge 1)\otimes E(\bar{\tau}_i \text{ }| \text{ }i\ge 0)$  which is defined on each algebra generator by the formulas
\[ 
	\begin{array}{c}
		\Delta (\bar{\xi}_n)=\underset{i+j=n}{\sum}\bar{\xi}_i\otimes\bar{\xi}_j^{p^i}\\
		\Delta(\bar{\tau}_n)=1\otimes \bar{\tau}_n + \underset{i+j=n}{\sum} \bar{\tau}_i \otimes \bar{\xi}_j^{p^i}.
	\end{array}
 \] 
 Here, by $\bar{x}$ we mean $\chi{x}$ where $\chi\co \mathcal{A}_*\rightarrow \mathcal{A}$ is the antipode structure map of the Hopf-algebra $\mathcal{A}_*$. When $\mathsf{E}=BP$, the co-product on elements of $BP_*BP\cong \mathbb{Z}_{(p)}[v_1,v_2, \dots]\otimes \mathbb{Z}_{(p)}[t_1,t_2, \dots]$  is defined by the formula
 \[ 
	\begin{array}{c}
		\Delta (t_n)=\underset{i+j=n}{\sum^F}t_i\otimes t_j^{p^i}
	\end{array}
 \] 
where $F$ is the formal group law of $BP$ associated to the  complex orientation $MU\rightarrow BP$, which equips $BP$ with the universal $p$-typical formal group law. 

We will write 
$X_p$ for the $p$-completion of a spectrum $X$, which agrees with the Bousfield localization $L_{S/p}X$ at the mod $p$ Moore spectrum $S/p$.  
We will write $\dot{=}$ to indicate that an equality holds up to multiplication by a unit in $\mathbb{F}_p$. Given a spectral sequence with signature
$E_n^{*,*}\implies A_*$,
we use the following terminology: 
\begin{enumerate}
\item An element $x\in E_n^{*,*}$ is an \emph{infinite cycle} if $d_r(x)=0$ for all $r\ge n$. 
\item An element $x\in E_n^{*,*}$ is a \emph{permanent cycle} if $d_r(x)=0$ for all $r\ge n$ and there does not exist a $y$ such that $d_r(y)=x$ for any $r\ge n$. 
\end{enumerate}
Also, throughout the paper, fix a prime $p\ge 5$ and let $q$ be a prime power that topologically generates $\mathbb{Z}_p^{\times}$, or equivalently $q$ generates the units in $\bZ/p^2\bZ$. 
\subsection{Acknowledgements}
This paper grew out of the author's Ph.D. thesis. The author would like to thank Andrew Salch for many discussions on the material in this paper and for his constant support and encouragement. Also, the author would like to thank Bob Bruner for offering his insight about the homological homotopy fixed point spectral sequence. The author also thanks W. Stephan Wilson for comments which led to a correction of an error in an early draft and Gerd Laures for feedback on the results of this paper in the context of the the double transfer and the $f$-invariant. 

\section{Overview of the toolkit} \label{prelim}
\subsection{The Hochschild--May spectral sequence for $\K(\mathbb{F}_q)_p$}
We recall necessary results and definitions from the author's paper \cite{K1localsphere} and the author's joint paper with A. Salch \cite{thhmay} since they will be cited later. We write for $\mathbb{N}^{\op}$ for the category with objects the natural numbers $n\ge 0$ and morphisms defined by 
\[ \Hom_{\mathbb{N}^{\text{op}}}(n,m)=\begin{cases} * & \text{ if }n\ge m \\ \emptyset & \text{ if }n <m. \end{cases}\]
\begin{defin} A filtered commutative ring spectrum $I$ is a cofibrant object in $\Comm \Sp^{\mathbb{N}^{\op}}$ where $\Comm \Sp^{\mathbb{N}^{\op}}$ has the model structure created by the forgetful functor to $\Sp^{\mathbb{N}^{\op}}$ and $\Sp^{\mathbb{N}^{\op}}$ has the projective model structure.\footnote{This definition differs slightly from that in \cite[Def 3.1.2]{thhmay}, but a cofibrant object in $\Comm \Sp^{\mathbb{N}^{\op}}$ is always a decreasingly filtered commutative monoid in $\Sp$ in the sense of loc. cit.}
 (See \cite[Sec. 4.1]{thhmay} for a discussion of why these model structures exist and have the desired properties). We write $I_i$ for $I$ evaluated on the natural number $i$. 
The associated graded of $I$ is a commutative ring spectrum $E_0^*I$ and it is defined so that, after forgetting the commutative monoid structure, it is the spectrum
$\bigvee_{i\ge 0} I_i/I_{i+1}$
where $I_i/I_{i+1}$ is the cofiber of the map $I_{i+1}\rightarrow I_i$. Note that by abuse of notation we may also regard the associated graded $E_0^*I=(I_0/I_1,I_1/I_2,\dots )$ as a commutative monoid in $\Sp^{\mathbb{N}^{\delta}}$ with respect to the Day convolution symmetric monoidal product where $\mathbb{N}^{\delta}$ is the natural numbers regarded as a discrete category. We can therefore consider topological Hochschild homology of $E_0^*I$ in the category $\Sp^{\mathbb{N}^{\delta}}$ so that $\THH(E_0^*I)$ is an object in $\Sp^{\mathbb{N}^{\delta}}$. Applying a generalized homology theory $\mathsf{E}$ level wise, we still have an object 
$\mathsf{E}\wedge \THH(E_0^*I)$ in $\Sp^{\mathbb{N}^{\delta}}$. Applying homotopy groups induces a functor
$ \pi_s^{\mathbb{N}^{\delta}} \colon \thinspace \Sp^{\mathbb{N}^{\delta}}\longrightarrow  \mathsf{Ab}^{\mathbb{N}^{\delta}}$
and then we write 
\begin{align}\label{notationexplanation} \mathsf{E}_{s,t}\THH(E_0^*I)=[\pi_s^{\mathbb{N}^{\delta}} (\mathsf{E}\wedge \THH(E_0^*I))](t).\end{align}
\end{defin}
\begin{exm}
By \cite[Thm 4.2.1]{thhmay}, an example of a filtered commutative ring spectrum associated to a connective commutative ring spectrum $R$ is the Whitehead filtration 
\[ \dots \rightarrow \tau_{\ge 2} R \rightarrow \tau_{\ge 1}R \rightarrow \tau_{\ge 0} R \]
which is equipped with structure maps $\rho_{i,j}:\tau_{\ge i} R\wedge \tau_{\ge j} R \rightarrow \tau_{\ge i+j}R$. Here $\tau_{\ge s}R$ is a spectrum with $\pi_i(\tau_{\ge s}R)\cong 0$ for $i<s$ that is equipped with a map $\tau_{\ge s}R\to R$ that induces an isomorphism on homotopy groups $\pi_i$ for $i\ge s$. 
We write simply $\tau_{\ge \bullet} R$ for the filtered commutative ring spectrum constructed in loc. cit. as a cofibrant object in $\Comm \Sp^{\mathbb{N}^{\op}}$. 
\end{exm} 
\begin{thm}[{\cite[Theorem 3.4.8 ]{thhmay}}] 
Let $I$ be a filtered commutative ring spectrum and let $\mathsf{E}$ be a connective spectrum. 
There is a spectral sequence with signature
\[ E^1_{s,t}=\mathsf{E}_{s,t}\THH(E_0^*I)\implies \mathsf{E}_{s}\THH(I_0), \] 
which we call the $\mathsf{E}$-Hochschild--May spectral sequence, where the bi-grading on the input is explained in Equation \eqref{notationexplanation}. 
\end{thm}
\begin{rem} 
When $I=\tau_{\ge \bullet} R$ we simply write $H\pi_*R$ for $E_0^*I$. It is a generalized Eilenberg-Maclane spectrum so whenever $\pi_kR$ is a finitely generated abelian group for all $k$ and 
 $\mathsf{E}=S/p$, $H\bF_p$, $V(1)$, or $BP\wedge V(1)$, then $\mathsf{E}_*\THH(H\pi_*R)$ is a graded $H\bF_p$-algebra and we can apply the following lemma to compute the input. 
\end{rem}
The following lemma is a consequence of the fact that all $H\bF_p$-modules are equivalent to a wedge of suspensions of $H\bF_p$ and an Adams spectral sequence argument (cf. \cite[Lem. 4.1]{MR2928844}).  
\begin{lem} \label{prim} Let $M$ be an $H\bF_p$-algebra. Then $M$ is equivalent to a wedge of suspensions of $H\bF_p$, and the Hurewicz map 
\[ \pi_*M\lra  H_*M\] 
induces an isomorphism between $\pi_*M$ and the subalgebra of $\mathcal{A}_*$-comodule primitives contained in $H_*M$. 
\end{lem} 
Using the lemma above, one can compute the $E^1$-page of the $H\mathbb{F}_p$-Hochschild--May spectral sequence, modulo an explicit description of $\HH_*(S/p_*(H\pi_*\K(\bF_q)_p))$ given in \cite[p. 271-272]{K1localsphere}. We will argue in the same way as \cite{K1localsphere}: there is an algebraic Hochschild--May spectral sequence which can be used to compute the $E^2$-page of the $\mathsf{E}$-Hochschild--May spectral sequence, for $\mathsf{E}=H\mathbb{F}_p$, $\mathsf{E}=H\mathbb{F}_p\wedge V(1)$ and $\mathsf{E}=V(1)$, and therefore at the $E^2$-page we have a simpler description. We give the description of both the $E^1$-page and the $E^2$-page. Note that the $E^1$-page of the $H\mathbb{F}_p$-Hochschild--May spectral sequence associated to the Whitehead filtration of $\K(\bF_q)_p$ is 
\[E_{*,*}^1=H_*(\THH(H\pi_*\K(\bF_q)_p)).\]
\begin{prop}[{\cite[Proposition 3.10, Proposition 3.11]{K1localsphere}}] \label{Bok1} 
There is an isomorphism of $\mathcal{A}_*$-comodule algebras
\begin{align*} 
E_{*,*}^1=H_*(\THH(H\pi_*\K(\bF_q)_p))\cong   
\cA//E(0)_*\otimes E(\sigma \bar{\xi}_1) \otimes P(\sigma \overline{\tau}_1) \otimes N 
\end{align*}
where $N=\HH_*(S/p_*(H\pi_*\K(\bF_q)_p))$ and 
where the $\cA_*$-co-action  on $\cA//E(0)_*$ is the usual one, that is the restriction of the co-product in $\mathcal{A}_*$, the coaction on $\sigma \bar{\xi}_1$ is 
\[ \psi(\sigma \bar{\xi}_1)=1\otimes \sigma \bar{\xi}_1\]
and the co-action on $\sigma \bar{\tau}_1$ is 
\[ \psi(\sigma \bar{\tau}_1)=1\otimes \sigma \bar{\tau}_1+\bar{\tau}_0\otimes \sigma \bar{\xi}_1.\]
The $E^2$-page of the $H\mathbb{F}_p$-Hochschild--May spectral sequence  associated to the Whitehead filtration of $\K(\bF_q)_p$ is 
\[ E_{*,*}^2= \mathcal{A}//E(0)_*\otimes E(\lambda_1)\otimes P(\mu_1)\otimes E(\alpha_1)\otimes P(v_1)\otimes E(\sigma v_1)\otimes \Gamma(\sigma \alpha_1)\]
as an $\mathcal{A}_*$-comodule, where $\lambda_1=\sigma \overline{\xi}_1$, $\mu_1=\sigma \overline{\tau}_1-\overline{\xi}_1\sigma \bar{\tau}_0$, with the usual $\mathcal{A}_*$-co-action on $\mathcal{A}//E(0)_*$ given by restricting the co-product on $\mathcal{A}_*$. The elements $\alpha_1v_1^k$, $v_1^k$, $\lambda_1$, $\mu_1$, and $\sigma \alpha_1$ are $\mathcal{A}_*$-comodule primitives. 
Multiplicatively, the $E^2$-page is isomorphic to 
\[ E_{*,*}^2= \mathcal{A}//E(0)_*\otimes  E(\lambda_1)\otimes P(\mu_1)\otimes  P(x_i:i\ge 1 )/(x_ix_j : i,j\ge 1)\otimes E(\sigma v_1)\otimes \Gamma(\sigma \alpha_1)\]
so that $x_{2i}=v_1^i$, $x_{2i-1}=\alpha_1v_1^{i-1}$ for $i\ge 1$.
The bidegrees of the indecomposable algebra generators at the $E^2$-page are given by 
 \begin{align*}
  |\sigma \bar{\xi}_1|&=(2p-1,0),\\ 
 |\sigma \bar{\tau}_1|&=(2p,0),\\ 
 |v_1^k|&=((2p-2)k,(2p-2)k-1),\\
 |\sigma v_1|&= (2p-1,2p-2),\\
 |\gamma_{k}(\sigma \alpha_1)|&=((2p-2)k,(2p-3)k),\\
  |\alpha_1v_1^{k-1}|&=((2p-2)k-1,(2p-2)k-1)\text{ for $k\ge 1$}. \footnotemark  
   \end{align*}
  The products $v_1^{k+1}\cdot v_1^{j+1}=0$, $\alpha_1 v_1^k \cdot v_1^{j+1}=0$, and $\alpha_1 v_1^k \cdot \alpha_1 v_1^j=0$ for $j,k\ge 0$ hold for bi-degree reasons as made clear by \cite[Figure 1]{K1localsphere}.
  \end{prop}
\footnotetext{There are a few minor typos in the statement of \cite[Prop. 3.11]{K1localsphere}: First, $\sigma \bar{\tau}_2$ should be $\sigma \bar{\tau}_1$. This is clear from \cite[Prop. 3.13]{K1localsphere}. Second, the ideal $\mathfrak{m}$ should be $(x_ix_j: i,j\ge 1)$. This is clear from the additive description. Third, we correct the bidegrees of the elements in the statement, where the correct bidegrees are clear from \cite[Figure 2]{K1localsphere}.} 
By \cite{MR0275429}, we know that 
\[ H_*(V(1))= E(\tau_0,\tau_1)\]
with $|\tau_0|=1$ and $|\tau_1|=2p-1$ with $\mathcal{A}_*$-coaction defined by 
\begin{align}
\label{coaction e0} \psi(\tau_0)=&1\otimes \tau_0+\tau_0\otimes 1\\
\label{coaction e1}  \psi(\tau_1)=&1\otimes \tau_1+\tau_1\otimes 1+\xi_1\otimes \tau_0.
 \end{align}
Note that using the recursive formulas 
\[\sum_{i=0}^n\xi_{n-i}^{p^{i}}\overline{\xi}_i=0\text{ and } \sum_{i=0}^n\xi_{n-i}^{p^{i}}\overline{\tau}_{i}+\tau_n=0\]
we can deduce that $\tau_0=-\overline{\tau}_0$, $\xi_1=-\bar{\xi}_1$, and consequently $\tau_1=-\overline{\tau}_1+\overline{\xi}_1\overline{\tau}_0$ (cf. \cite[Lem. 10]{MR0099653}). This allows us to rewrite the formulas above in terms of our chosen basis for $\mathcal{A}_*$. 

Combining this with Proposition \ref{Bok1}, we computed the $E^1$-page and $E^2$-page of  the $V(1)$-Hochschild--May spectral sequence using Lemma \ref{prim}, which we again recall modulo  an explicit description of $\HH_*(S/p_*(H\pi_*\K(\bF_q)_p))$, which can be found in \cite[p. 271-272]{K1localsphere}. Recall that the $E^1$-page of the $V(1)$-Hochschild--May spectral sequence associated to the Whitehead filtration of $\K(\bF_q)_p$ is 
\[E^1_{*,*}=V(1)_*\THH(H\pi_*\K(\bF_q)_p).\]
\begin{prop}[{\cite[Prop. 3.15]{K1localsphere}}]\label{computation}
There is an isomorphism of graded $\mathbb{F}_p$-algebras
 \[ E_{*,*}^1=V(1)_*\THH(H\pi_*\K(\bF_q)_p)\cong  E(\lambda_1, \epsilon_1)\otimes P(\mu_1)\otimes \HH_*(S/p_*(H\pi_*\K(\bF_q)_p))\] 
and the $E^2$-page of the $V(1)$-Hochschild--May spectral sequence associated to the Whitehead filtration of $\K(\bF_q)_p$ is 
 \[ E^2_{*,*}=E(\lambda_1,\epsilon_1)\otimes P(\mu_1)\otimes E(\alpha_1)\otimes P(v_1)\otimes \Gamma(\sigma \alpha_1)\otimes E(\sigma v_1)\]
 as a graded $\mathbb{F}_p$-module where 
 \begin{align*} 
 |\epsilon_1|= |\lambda_1|&=(2p-1,0),\\ 
 |\mu_1|&=(2p,0),\\ 
 |v_1^k|&=((2p-2)k,(2p-2)k-1),\\
 |\sigma v_1|&= (2p-1,2p-2),\\
 |\gamma_k(\sigma \alpha_1)|&=((2p-2)k,(2p-3)k),\\
  |\alpha_1v_1^{k-1}|&=((2p-2)k-1,(2p-2)k-1)  \text{ for $k\ge 1$},
 \end{align*}
Here $\epsilon_1=\tau_1-\tilde{\tau}_1$ where  $\tilde{\tau}_1$ is the indecomposable algebra generator in degree $2p-1$ in $\mathcal{A}//E(0)_*\subset H_*(\THH(H\pi_*\K(\bF_q)_p))$ and the multiplicative structure on the $E^2$-page of the $V(1)$-Hochschild--May spectral sequence associated to the Whitehead filtration of $\K(\bF_q)_p$ is
 \[ E^2_{*,*}=E(\lambda_1,\epsilon_1)\otimes P(\mu_1)\otimes \Gamma(\sigma \alpha_1)\otimes E(\sigma v_1)\otimes P(x_i: i\ge 1)/(x_ix_j:i,j\ge 1)\]
where $x_{2i}=v_1^{i}$ and $x_{2i-1}=\alpha_1v_1^{i-1}$ for $i\ge 1$. 
\end{prop}

Our computations build on the computation of homology of topological Hochschild homology of $\K(\bF_q)_p$ due to  Angeltveit--Rognes \cite{MR2171809}. 
\begin{thm}[Theorem 7.13 and Theorem 7.15 \cite{MR2171809}] \label{HFpj}
There is an isomorphism of $\mathcal{A}_*$-comodule algebras
\[ H_*(\K(\bF_q)_p)\cong P(\tilde{\xi}_1^p,\tilde{\xi}_2,\bar{\xi}_3,...)\otimes E(\tilde{\tau}_2,\bar{\tau}_3, ...)\otimes E(b) \cong (\mathcal{A}//A(1))_*\otimes E(b)\] 
where all the elements in $(\mathcal{A}//A(1))_*$ besides $\tilde{\tau}_2$, $\tilde{\xi}_1^p$, and $\tilde{\xi}_2$, and $b$ have the usual $\mathcal{A}_*$-co-action and the co-action on the remaining elements $\tilde{\tau}_2$, $\tilde{\xi}_1^p$, $\tilde{\xi}_2$, and $b$ are
\begin{align*} 
\psi(b)=&1\otimes b \\
\psi(\tilde{\xi}_1^p)=& 1\otimes \tilde{\xi}_1^p -\tau_0\otimes b + \bar{\xi}_1^p\otimes 1 \\
\psi(\tilde{\xi}_2)= &1\otimes \tilde{\xi}_2+\bar{\xi}_1\otimes \tilde{\xi}_1^p +\tau_1\otimes b +  \bar{\xi}_2\otimes 1 \\
\psi(\tilde{\tau}_2)=& 1\otimes \tilde{\tau}_2 +\bar{\tau}_1\otimes \tilde{\xi}_1^p + \bar{\tau}_0\otimes \tilde{\xi}_2 - \tau_1 \tau_0\otimes b  + \bar{\tau}_2\otimes 1  .\\
\end{align*}

The degrees of $\tilde{\tau}_2$, $\tilde{\xi}_1^p$, and $\tilde{\xi}_2$, and $b$ are 
\begin{align*}
	|\tilde{\tau}_2|&=2p^2-1,\\ 
	|\tilde{\xi}_1^p|&=2p^2-2p, \\
	|\tilde{\xi}_2|&=2p^2-2, \text{ and }\\  
	|b|&=2p^2-2p-1
\end{align*}
respectively. There is also an isomorphism
\[ H_*(\THH(\K(\bF_q)_p))\cong H_*(\K(\bF_q)_p)\otimes E(\sigma \tilde{\xi}_1^p,\sigma \tilde{\xi}_2)\otimes P(\sigma \tilde{\tau}_2)\otimes \Gamma(\sigma b) \]
of $\mathcal{A}_*$-comodules and $H_*(\K(\bF_q)_p)$-algebras. The $\mathcal{A}_*$-co-action is given by using the formula
\[ \psi(\sigma x)=(1\otimes \sigma)\circ \psi(x) \]
and the previously stated co-actions. 
\end{thm}
  \begin{cor}\label{homologymodpv1}
 There is an isomorphism of $\mathcal{A}_*$-comodule algebras
 \[ H_*(V(1)\wedge \THH(\K(\bF_q)_p))\cong E(\tau_0,\tau_1)\otimes H_*(\THH(\K(\bF_q)_p))\]
 with $\mathcal{A}_*$-co-action on $\tau_0$ given in \eqref{coaction e0}, the coaction on $\tau_1$ given in \eqref{coaction e1}, and the coaction on $x\in H_*(\THH(\K(\bF_q)_p))$ given in Theorem \ref{HFpj}. 
 \end{cor}

Finally, we recall the computation of $V(1)$-homotopy of topological Hochschild homology of $\K(\bF_q)_p$. 
This was computed by the author in \cite[Thm. 1.5]{K1localsphere} and by H\"oning in \cite[Thm. 7.2]{Hon19} using different methods. \footnotemark 
\begin{thm} \label{mod p v_1 THH K}
There is an isomorphism of graded $\mathbb{F}_p$-algebras
\[ V(1)_*\THH(\K(\bF_q)_p)\cong P(\mu_2)\otimes \Gamma(\sigma b)\otimes\bF_p\{1,\alpha_1,\lambda_1',\lambda_2\alpha_1, \lambda_2\lambda_1', \lambda_2\lambda_1'\alpha_1\}.\] 
where $\alpha_1\cdot  (\lambda_2\lambda_1')=\lambda_1' \cdot (\lambda_2\alpha_1)= \lambda_2\lambda_1'\alpha_1$.
\end{thm}
\footnotetext{There is a minor typo in the statement of \cite[Thm. 1.5]{K1localsphere}. The element $1$ should also be included. This is evident since $V(1)_*\THH(\K(\bF_q)_p)$ is ring.}

\subsection{The homotopy fixed point spectral sequence}
In this section, we summarize results from \cite[\S 2-4]{BR05}, which also apply when replacing $H\mathbb{F}_p$ with any spectrum $\mathsf{E}$. 
\begin{defin}
We write 
\[\mathsf{E}_{s+t}^c(X^{h\bT}) := \lim \mathsf{E}_{s+t}F(S(\bC^k)_+,X)^\bT. \]
for the \emph{continuous $\mathsf{E}$-homology of $X^{h\bT}$} and we refer to the filtration $F(S(\bC^k)_+,X)^\bT$ as the \emph{skeletal} filtration. 
When $\mathsf{E}$ is a finite spectrum, then $\mathsf{E}_{*}^c(X^{h\bT})=\mathsf{E}_{*}(X^{h\bT})$ since smashing with a dualizable spectrum commutes with homotopy limits. 
\end{defin}

\begin{prop}[c.f. Proposition 2.1 and Proposition 4.1 in \cite{BR05}] \label{BRSS}
Let $\mathsf{E}$ be a homotopy associative ring spectrum. There is a spectral sequence $\mathsf{E}_*\mathsf{E}$ comodules with signature 
\[ E^2_{*,*}=P_{\mathbb{Z}}(t)\otimes \mathsf{E}_*X \implies \mathsf{E}_*^c(X^{h\bT}),\]
which we simply call the \emph{homotopy fixed point spectral sequence}. If, in addition, $X$ is a commutative ring spectrum, then this is a spectral sequence of $\mathsf{E}_*\mathsf{E}$-comodule algebras where $\mathsf{E}_*X$ has the Pontryagin product. 
\end{prop}
\begin{proof} 
The proof is essentially the same as in \cite[Prop. 2.1, Prop. 4.1]{BR05} so we omit it. 
\end{proof}
\begin{lem}\label{prop:d2homofixedpoint}
The $d^2$ differentials in the homotopy fixed point spectral sequence associated to $\THH(R)$ are of the form 
\begin{align*} d^2(x)=t\sigma x. \end{align*}
where $t$ is the generator of $H_{gp}^{-*}(\bT;\mathbb{Z})\cong P_{\mathbb{Z}}(t)$ in degree $-2$ whenever $\eta\in \pi_1S$ acts trivially on $\mathsf{E}_*\THH(R)$. 
\end{lem}
\begin{proof}
This follows directly from \cite[Lemma 1.4.2]{Hes96}. 
\end{proof}
\begin{defin}
We let $T_k(R)$ denote the spectrum  $F(S(\mathbb{C}^{k+1})_+,\THH(R))^{\bT},$
which has the property that $\holim_k T_k(R)=\THH(R)^{h\bT}$. We write 
\[ \TC^{-}(R):=\THH(R)^{h\bT}\]
following the convention of \cite{2016arXiv160201980H} and we also write 
\[\mathsf{E}_*^c(\TC^{-}(R)):=\mathsf{E}_*^c(\THH(R)^{h\bT}).\]
Note that $T_0(R)=\THH(R)$. 
\end{defin}
\begin{remark}\label{truncated ss}
Note that there is a \emph{truncated homotopy fixed point spectral sequence} with signature 
\[ E^2_{*,*}=P_{\mathbb{Z}}(t)/t^{k+1}\otimes \mathsf{E}_*\THH(R) \implies \mathsf{E}_*T_k(R)\] 
that is constructed in exactly the same way as the $\mathbb{T}$-homotopy fixed point spectral sequence. The $d^2$ differential is the same as \ref{prop:d2homofixedpoint} whenever $\eta\in \pi_1S$ acts trivially on $\mathsf{E}_*\THH(R)$. 
\end{remark}

\section{Detecting $\beta$ elements in iterated algebraic K-theory of finite fields}
\subsection{The $\mathsf{E}$-Adams spectral sequence}\label{ASS}
The following discussion applies to any homotopy commutative and homotopy associative ring spectrum $\mathsf{E}$ such that $(\mathsf{E}_*,\mathsf{E}_*\mathsf{E})$ is a flat Hopf algebroid satisfying some mild additional hypotheses \cite[Assumption 2.2.5]{rav1}, but for simplicity we assume $\mathsf{E}=H\mathbb{F}_p$ or $\mathsf{E}=BP$ (cf. \cite[Prop. 2.2.6]{rav1}). We write $\Delta \colon \thinspace \mathsf{E}_*\mathsf{E}\to \mathsf{E}_*\mathsf{E}\otimes_{\mathsf{E}_*}\mathsf{E}_*\mathsf{E}$ for the co-product and $\eta_R\colon \thinspace \mathsf{E}_*\to \mathsf{E}_*\mathsf{E}$ for the right unit. Given a spectrum $X$ we write $\psi \colon \thinspace \mathsf{E}_*X\to \mathsf{E}_*\mathsf{E}\otimes_{\mathsf{E}_*}\mathsf{E}_*X$ for the left $(\mathsf{E}_*,\mathsf{E}_*\mathsf{E})$-coaction on $\mathsf{E}_*X$. We also assume $X$ is bounded below and $\mathsf{E}_*X$ is finite type for simplicity. 

We briefly recall the construction of the $\mathsf{E}$-Adams spectral sequence. For more details, we refer the reader to \cite[Ch. 2 \S 2]{rav1}. The main examples we will consider are $\mathsf{E}=H\mathbb{F}_p$ and $\mathsf{E}=BP$. Classically, we can build the $\mathsf{E}$-Adams spectral sequence for $X$ using the canonical $\mathsf{E}$-Adams resolution 
\begin{align}\label{Adams tower}
	\xymatrix{ 
			X \ar[d]_{f_0} &  \ar[l]_{g_0}  \Sigma^{-1}\overline{\mathsf{E}}\wedge X  \ar[d]_{f_1} &  \ar[l]_{g_1}\Sigma^{-2}\overline{\mathsf{E}}^{\wedge 2} \wedge X  \ar[d]_{f_2} &\ar[l]_(.4){g_2 } \dots \\
			\mathsf{E}\wedge X \ar@{-->}[ur]^{h_0}&  \Sigma^{-1}\mathsf{E} \wedge \overline{\mathsf{E}}\wedge X \ar@{-->}[ur]^{h_1}& \Sigma^{-2}\mathsf{E} \wedge \overline{\mathsf{E}}^{\wedge 2}\wedge X \ar@{-->}[ur]^{h_2}&&
	}
\end{align}
from \cite[Lem. 2.2.9]{rav1} where $\overline{\mathsf{E}}$ is the cofiber of the unit map $\eta_{\mathsf{E}} \colon \thinspace S\to \mathsf{E}$, $f_i=\eta_\mathsf{E}\wedge \Sigma^{-i} \overline{\mathsf{E}}^{\wedge i}\wedge X$ and $g_i=(\Sigma^{-1}\overline{\mathsf{E}}\to S)\wedge \Sigma^{-i}\overline{\mathsf{E}}^{\wedge i}\wedge X$ and $h_i\colon \thinspace \mathsf{E}\wedge X\to \overline{\mathsf{E}}\wedge X$ induces the boundary map in the long exact sequence in homotopy groups. 
Applying $\pi_*$ produces an exact couple with $E_1$-page given by the differential graded algebra 
\[ \mathsf{E}_*X\longrightarrow   \mathsf{E}_*\overline{\mathsf{E}}\otimes_{\mathsf{E}_*}  \mathsf{E}_*X\longrightarrow \mathsf{E}_* \overline{\mathsf{E}}^{\otimes_{\mathsf{E}_*}2}\otimes_{\mathsf{E}_*}  \mathsf{E}_*X\longrightarrow \dots \]
called the (normalized) \emph{cobar complex} with boundary maps induced by the maps $\Sigma^{-i}h_i\circ \Sigma^if_i$. Specifically, the differential 
\[ d_1^{s,t}\colon \thinspace [  \mathsf{E}_* \overline{\mathsf{E}}^{\otimes_{\mathsf{E}_*}s}\otimes_{\mathsf{E}_*} \mathsf{E}_*X]_t\to  [ \mathsf{E}_* \overline{\mathsf{E}}^{\otimes_{\mathsf{E}_*} s+1}\otimes_{\mathsf{E}_*} \mathsf{E}_*X]_t\]
in the cobar complex is given by the formula
\begin{align*}
d_1^{s,t}(e_1\otimes \dots \otimes e_s\otimes m)=& \sum_{i=1}^{s} (-1)^{i} e_1\otimes \dots \otimes e_{i-1}\otimes \overline{\Delta}(e_i)\otimes e_{i+1}\otimes \dots \otimes e_s\otimes m\\
 & +(-1)^{s+1}e_1\otimes \dots \otimes e_s\otimes \overline{\psi}(m)
\end{align*} 
when $e_1\otimes \dots \otimes e_s\otimes m\in  [\mathsf{E}_* \overline{\mathsf{E}}^{\otimes_{\mathsf{E}_*} s}\otimes_{\mathsf{E}_*} \mathsf{E}_*X]_t$, where $\overline{\Delta}(x)=\Delta(x)-x\otimes 1-1\otimes x$ and $\overline{\psi}(m)=\psi(m)- 1\otimes m$. 
We note that under our hypotheses on $\mathsf{E}$ and $X$ the $\mathsf{E}$-Adams spectral sequence has signature
\[E_2^{s,t}= [H^s( \mathsf{E}_* \overline{\mathsf{E}}^{\otimes_{\mathsf{E}_*} \bullet}\otimes_{\mathsf{E}_*}  \mathsf{E}_*X; d_1)]_t\implies \pi_{t-s}X_{p}\]
and it strongly converges. 
\begin{remark}
If we consider the $\mathsf{E}$-Adams spectral sequence when $X$ is the sphere spectrum we drop the last tensor factor from the cobar complex. 
\end{remark}

\subsection{Detecting $v_2$ and $\beta_1$}
The mod $p$ Moore spectrum $S/p$ and the Smith-Toda complex $V(1)$ are defined so that they fit into exact triangles 
\begin{align}\label{exact triangle 1} \xymatrix{ S \ar[r]^p & S \ar[r]^{i_0} & S/p \ar[r]^{j_0} & \Sigma S }\end{align}
and 
\begin{align}\label{exact triangle 2} \xymatrix{ \Sigma^{2p-2}S/p \ar[r]^{v_1} & S/p \ar[r]^{i_1} & V(1) \ar[r]^{j_1} & \Sigma^{2p-1}S/p }\end{align}
in the stable homotopy category of spectra. We will abuse notation and write  $i_0\co \pi_*S\to S/p$ and $i_1\co \pi_*S/p\to \pi_*V(1)$ for the maps induced by $i_0$ and $i_1$ respectively. Similarly, we write $j_0 \colon \thinspace \pi_{*}S/p\to \pi_*\Sigma S$ and $j_1\colon \thinspace \pi_*V(1)\to \pi_*\Sigma^{2p-1}S/p$ for the maps induced by $j_0$ and $j_1$. We similarly abuse notation for composites of these maps. Throughout this section we write $\overline{\mathcal{A}}=\pi_*(H\mathbb{F}_p\wedge \overline{H\mathbb{F}}_p)$ for the kernel of the augmentation $\mathcal{A}\longrightarrow \mathbb{F}_p$. 
\begin{defin}\label{elements}
Let 
\begin{align*}
v_2\in \pi_{2p^2-2}V(1)  \\ 
\beta_1^{\prime}\in \pi_{2p^2-2p-1}V(1) \\
\beta_1\in \pi_{2p^2-2p-2}V(1) &
\end{align*}
be the classes represented in the respective cobar complexes by the permanent cycles
\begin{align*}
	\overline{v}_2:=&\tau_2\otimes 1+\xi_2\otimes \tau_0+\xi_1^p\otimes \tau_1 \in \overline{\mathcal{A}}_*\otimes H_*(V(1)), \\
	h_1:=&\bar{\xi}_1^p \otimes 1\in \overline{\mathcal{A}}_*\otimes H_*(V(1)), \text{ and }\\
	b_0:=&\sum_{j=1}^{p-1}\frac{1}{p}\binom{p}{j} \bar{\xi}_1^{p-j}\otimes \bar{\xi}_1^{j}\otimes 1  \in  \overline{\mathcal{A}}^{\otimes 2}_*\otimes H_*(V(1)).
\end{align*}
The bidegrees of these elements are 
\begin{align*}
|\overline{v}_2|=&(1,2p^2-1),\\ 
|h_1|=& (1,2p^2-2p), \text{ and }\\
|b_0|=&(2,2p^2-2p)
\end{align*}
where $|x|=(s,t)$ corresponds to stem $t-s$ and Adams filtration $s$. 
\end{defin}
The following lemma is immediate from Proposition \ref{prop:d2homofixedpoint}:
\begin{lem}\label{little lem}
There is an isomorphism 
\[ H_*(T_1(\K(\mathbb{F}_q)_p))\cong H_*(H_*(\THH(\K(\mathbb{F}_q)_p)\otimes E(t);d)\]
where the right-hand side is the homology of the differential graded algebra $H_*(\THH(\K(\mathbb{F}_q)_p))\otimes E(t)$ with $|t|=(0,-2)$ and $|x|=(s,0)$ if $x\in H_s(\THH(\K(\mathbb{F}_q)_p))$. The differential is $d^1(x)=t\sigma x$ for all $x\in H_*(\THH(\K(\mathbb{F}_q)_p)$.
\end{lem}
In the proof of the following proposition, we will will make use of differentials in both the homotopy fixed point spectral sequence and the Adams spectral sequence. To differentiate between the two, we use notation $d^r$ for differentials in the homotopy fixed point spectral sequence and we use the notation $d_r$ for differentials in the Adams spectral sequence. 
The following argument is inspired by an argument of Ausoni--Rognes \cite[Prop. 4.8]{MR1947457}. 
\begin{prop} \label{perm classes} 
The classes $v_2$, $i_1i_0\beta_1$, and $i_1\beta_1^{\prime}$ in $V(1)_*$ map nontrivially to the classes $t\mu_2$, $t\sigma b$, and $t\lambda_1^{\prime}$ respectively in 
$V(1)_{*}\TC^{-}(\K(\bF_q)_p)$ up to multiplication by a unit in $\mathbb{F}_p$. 
\end{prop} 

\begin{proof}
We first sketch the strategy. First of all, it will suffice to prove that each of the classes $v_2$, $i_1i_0\beta_1$, and $i_1\beta_1^{\prime}$ map to non-trivial elements $t\mu_2$, $t\sigma b$, and $t\lambda_1^{\prime}$ in the homotopy groups of the skeleton $V(1)\wedge T_1(\K(\bF_q)_p$ since the unit map $V(1)\to V(1)\wedge T_1(\K(\bF_q)_p$ factors through the unit map $V(1)\to V(1)\wedge \TC^{-}(\K(\bF_q)_p)$. We compute $V(1)_*T_1(\K(\bF_q)_p$ using two different spectral sequences: the $H\mathbb{F}_p$-Adams spectral sequence from Section \ref{ASS} and the truncated homotopy fixed point spectral sequence from Remark \ref{truncated ss}. By representing each of the classes $v_2$, $i_1i_0\beta_1$, and $i_1\beta_1^{\prime}$ by explicit classes in the $E_1$-page of the Adams spectral sequence, we compute the image in the $E_1$-page of the Adams spectral sequence for $V(1)\wedge T_1(\K(\bF_q)_p$ and we compute directly that these classes are permanent cycles. Using the truncated homotopy fixed point spectral sequence from  Remark \ref{truncated ss}, we compute that there is only one class in the correct degree in $V(1)_*T_1(\K(\bF_q)_p$ for each of the classes $v_2$, $i_1i_0\beta_1$, and $i_1\beta_1^{\prime}$. This determines that the explicit names for the images of each of the classes $v_2$, $i_1i_0\beta_1$, and $i_1\beta_1^{\prime}$ must be $t\mu_2$, $t\sigma b$, and $t\lambda_1^{\prime}$ up to multiplication by a unit in $\mathbb{F}_p$. 

First, $v_2$ is represented by $\bar{v}_2$, $\beta_1'$ is represented by $h_1$ and $\beta_1$ is represented by $b_{0}$
in  the $E_1$-page of the Adams spectral sequence that converges to $\pi_*V(1)$ as in Definition \ref{elements}. We consider the map of Adams spectral sequences 
\[ \Ext_{\mathcal{A}_*}^{*,*}(\bF_p, H_*(V(1)))\lra \Ext_{\mathcal{A}_*}^{*,*}(\bF_p, H_*(V(1))\otimes H_*( T_1(\K(\bF_q)_p) ) \] 
induced by the unit map
\[ \xymatrix{ V(1)\wedge S \ar[rr]^(.4){1_{V(1)}\wedge \eta }&& V(1)\wedge T_1(\K(\bF_q)_p)   )}.\]
We see that $\overline{v}_2$, $h_1$, and $b_{0}$ are permanent cycles in the source,
which map to classes of the same name in the target. Since the elements in the source are infinite cycles, this implies that the elements that they map to are infinite cycles as well.  We therefore just have to check that these classes are not boundaries. 

We can eliminate the possibility of a $d_1$ differential with $\overline{v}_2$ as a boundary by computing the differential in the cobar complex for $H_*(V(1)\wedge T_1(\K(\bF_q)_p))$ on each class of the correct degree. There is a truncated homotopy fixed point spectral sequence by Remark \ref{truncated ss} with signature
\[ E^2_{*,*}= E(t)\otimes H_*(V(1)\wedge THH(\K(\bF_q)_p))\implies H_*(V(1))\otimes H_*(T_1(\K(\bF_q)_p))\]
where  
\[E^2_{*,*}=E(t)\otimes E(\tau_0,\tau_1)\otimes P(\tilde{\xi}_1^p,\tilde{\xi}_2,\overline{\xi}_3,\dots )\otimes E(b,\tilde{\tau}_2,\overline{\tau}_3,\dots)\otimes E(\sigma \tilde{\xi}_1^p, \sigma \tilde{\xi}_2)\otimes P(\sigma \tilde{\tau}_2)\otimes \Gamma(\sigma b)\]
by Corollary \ref{homologymodpv1}. There is an isomorphism 
\[ E^4_{*,*}=H_*(V(1)) \otimes H_*(H_*(\THH(\K(\mathbb{F}_q)_p)\otimes E(t);d^2)\]
where $d^2(x)=t\sigma x$ and the spectral sequence clearly collapses at the $E^4$-page for bi-degree reasons as described in Lemma \ref{little lem}. Consequently, 
\[ H_{2p^2-1}(V(1)\wedge T_1(\K(\bF_q)_p))\cong \mathbb{F}_p\{\sigma \tilde{\xi}_2,\tau_1\sigma b\}\]
Since $\sigma b$ is primitive, $d_1(\sigma b)=0$. The coaction on  $\tau_1$ is given by the co-product 
\[ \psi(\tau_1)=\tau_1\otimes 1 + \xi_1\otimes \tau_0 +1\otimes \tau_1\]
and the co-action on $\sigma \tilde{\xi}_2$ is described in \ref{HFpj}. 
Consequently, $d_1(\tau_1)=\tau_1\otimes 1+\xi_1\otimes \tau_0$ and we compute 
\begin{align*}
d_1(\sigma \tilde{\xi}_2)=& \bar{\xi}_1\otimes \sigma \tilde{\xi}_1^p +\tau_1\otimes \sigma b, \text{ and }\\
d_1(\tau_1\sigma b)=& \tau_1\otimes \sigma b+ \xi_1 \otimes \tau_0 \sigma b
\end{align*} 
and then observe that 
\begin{align*}
	d_1(\alpha \tau_1\sigma b+\beta \sigma \tilde{\xi}_2)=& \alpha (\tau_1\otimes \sigma b+\xi_1 \otimes  \tau_0 \sigma b)+\beta(\bar{\xi}_1\otimes \sigma \tilde{\xi}_1^p +\overline{\tau}_1\otimes \sigma b)\\
	\ne & \tau_2\otimes 1+\xi_2\otimes \tau_0+\xi_1^p\otimes \tau_1
\end{align*}
for any $\alpha,\beta\in \mathbb{F}_p$.  
Therefore, $\overline{v}_2$ survives to the $E_{2}$-page. 
There are no possible longer differentials hitting $\overline{v}_2$ because $\overline{v}_2$ is in Adams filtration one; hence, it is a permanent cycle. 

We eliminate the possibility that the class $h_1$ is a boundary of a $d_1$ by the same method. The source of a $d_1$-differential hitting $h_1$ must be in Adams filtration zero and degree $2p^2-2p$. 
By Proposition \ref{prop:d2homofixedpoint}, we compute $H_{2p^2-2p}(V(1)\wedge T_1(\K(\bF_q)_p))\cong \mathbb{F}_p\{\sigma b\}$
we just need to compute $d_1(\sigma b)$. However, we already observed that $d_1(\sigma b)=0$. 
The class $h_1$ is in Adams filtration one so it cannot be the target of a longer differential, therefore it is a permanent cycle. 

For $b_{0}$, we need to check that it is not the boundary of a $d_1$ or a $d_2$, because it is in Adams filtration two. We first need to check that  $b_{0}$, is not a $d_1$-boundary for any element in 
$\overline{\mathcal{A}}_*\otimes H_*(V(1)\wedge T_1(\K(\bF_q)_p))$ 
in degree $2p^2-2p$. 
We check the differential in the cobar complex on all the elements here in the right degree. These classes are 
\[
 	\left \{ 
		\begin{array}{cccc}
			\bar{\tau}_0 \otimes \tau_0 t\tilde{\xi}_1^p, &
			\bar{\xi}_1^{p-1}\otimes \tau_0 \tau_1t, &
			\bar{\xi}_1^{p-2}\bar{\tau}_1 \otimes \tau_1t, &
			\bar{\xi}_1^{p-1}\bar{\tau}_0\otimes \tau_1t,\\
			\bar{\xi}_1^{p-1}\bar{\tau}_1\otimes \tau_0t, &
			\bar{\xi}_1^{p-1}\bar{\tau}_0\bar{\tau}_1\otimes t,&
			\bar{\xi}_1^p \otimes 1 &
		\end{array}
 	\right \} .
 \] 

Recall that Milnor computed the co-action of $\mathcal{A}_*$ on 
\[H^*(\mathbb{C}P^{\infty},\bF_p)=H^*(B\bT; \bF_p) =H_{\text{gp}}^*(\bT;\bF_p),\]
and the co-action on the class $t$ is 
\[\psi(t)=\sum_{i\ge 0} \bar{\xi}_i\otimes t^{p^i},\]
see \cite{MR0099653}. 
Therefore, in the input of the truncated homotopy fixed point spectral sequence computing $V(1)_*T_1(\K(\bF_q)_p)$, the $\mathcal{A}_*$ co-action on $t$ is primitive. 

We compute the differential in the cobar complex on each of the elements that could possibly have the class representing $\beta_1$ as a target. However, since all of the classes except $\bar{\xi}_1^p \otimes 1$ have a $t$ as a factor of the rightmost tensor factor and $t$ is a comodule primitive, it suffices to compute the $d_1$-differential on the class $\bar{\xi}_1^p \otimes 1$. We then compute directly that 
\[ d_1(\bar{\xi}_1^p \otimes 1 ) = 0.\]

We can conclude that $b_{0}$ is not a boundary of a $d_1$. 

Since $b_{0}$ is in Adams filtration two, we still have to check that there is no $d_2$ differential hitting it in the Adams spectral sequence, 
\[ \Ext_{\mathcal{A}_*}^{*,*} (\bF_p , H_*(V(1)\wedge T_1(\K(\bF_q)_p) ) \implies V(1)_*T_1(\K(\bF_q)_p). \] 
Since a $d_2$ would have to have its source on the $0$-line in degree $2p^2-2p-1$, it would have to be a class in $H_{2p^2-2p-1}(V(1)\wedge T_1(\K(\bF_q)_p))$. 

Again, we compute 
\[ H_{2p^2-2p-1}(V(1)\wedge T_1(\K(\bF_q)_p))\cong \bF_p\{\tau_0t\tilde{\xi}_1^p\}\] 
by Proposition \ref{prop:d2homofixedpoint}. 
Since $d_1(\tau_0)=-\tau_0\otimes 1$, the Leibniz rule implies 
\[ d_1(\tau_0t\tilde{\xi}_1^p)=-\tau_0\otimes t\tilde{\xi}_1^p - \tau_0 d_1(t\tilde{\xi}_1^p),\] 
so since $d_1(t\tilde{\xi}_1^p)=-\bar{\xi}_1^p\otimes t-\bar{\tau}_0\otimes tb$, we observe that 
\[d_1(\tau_0t\tilde{\xi}_1^p)=-\tau_0\otimes  t\tilde{\xi}_1^p +  \tau_0 \bar{\xi}_1^p\otimes t\]  
So $\tau_0(t\tilde{\xi}_1^p)$ does not survive to the $E_2$-page and therefore it cannot support a $d_2$ differential hitting $b_{0}$.  
Therefore, since $b_{0}$ is an infinite cycle and we have proven it isn't a boundary
it must be a permanent cycle.  

We conclude the elements $v_2$, $\beta_1$ and $\beta_1^{\prime}$ 
map nontrivially from 
$ V(1)_*$ to $V(1)_* T_1(\K(\bF_q)_p)$
via the map induced by the unit map $S\to T_1(\K(\bF_q)_p)$. We claim that in $V(1)_*T_1(\K(\bF_q)_p)$, the only possible classes in the correct degree to be $v_2$, $\beta_1$ and $\beta_1^{\prime}$ are  $t\mu_2$, $t\sigma b$ and $t\lambda_1^{\prime}$, respectively. To see this consider the truncated homotopy fixed point spectral sequence 
\[ E(t)\otimes V(1)_*THH(\K(\bF_q)_p) \implies V(1)_*T_1(\K(\bF_q)_p)\]
the only possible $d^2$-differential in this spectral sequence  that could affect degree $2p^2-2p-1$, $2p^2-2p-2$, or $2p^2-2$ is $d^2(\lambda_1^{\prime})=a t\sigma b$ for some $a\in \mathbb{F}_p$. We can show that $a=0$ by considering the map of truncated homotopy fixed point spectral sequences with codomain
\[ E(t)\otimes H_*(V(1)\wedge THH(\K(\bF_q)_p) \implies H_*(V(1)\wedge T_1(\K(\bF_q)_p)\]
which sends $\lambda_1^{\prime}$ and $t\sigma b$ to classes with the same name. We then observe that $d^2(\lambda_1^{\prime})=0$ in the codomain spectral sequence. Consequently, we have isomorphisms 
\[ V(1)_{2p^2-2}T_1(\K(\bF_q)_p)\cong \mathbb{F}_p\{t\mu_2\}\]
\[ V(1)_{2p^2-2p-2}T_1(\K(\bF_q)_p)\cong \mathbb{F}_p\{t\sigma b\}\]
\[ V(1)_{2p^2-2p-1}T_1(\K(\bF_q)_p)\cong \mathbb{F}_p\{t\lambda_1^{\prime}\}\]
proving the claim. 
The unit map $S\to T_1(\K(\bF_q)_p)$ factors through the unit map $S\to \TC^{-}(\K(\bF_q)_p)$, so these classes pull back to classes in $V(1)_*\TC^{-}(\K(\bF_q)_p)$. 
\end{proof}

\begin{cor}\label{perm cycles homotopy fixed}
The classes $t\mu_2$, $t\sigma b$, and $t\lambda_1^{\prime}$ are permanent cycles in the homotopy fixed point spectral sequence 
\[ H^*(\bT,V(1)_*\THH(\K(\mathbb{F}_q)_p))\implies V(1)_*\TC^{-}(\K(\mathbb{F}_q)).\]
\end{cor}
\begin{lem}\label{d2p-2}
There is a differential $d_{2p-2}(t)=t^p\alpha_1$ in the homotopy fixed point spectral sequence 
\[ H^*(\bT,V(1)_*\THH(\K(\mathbb{F}_q)_p))\implies V(1)_*\TC^{-}(\K(\mathbb{F}_q)).\]
\end{lem}
\begin{proof}
First, we can show that there is a differential $d_{2p-2}(t)=t^p\alpha_1$ in the homotopy fixed point spectral sequence 
\[ H^*(\bT,V(1)_*\K(\mathbb{F}_q)_p))\implies V(1)_*\K(\mathbb{F}_q)_p^{h\mathbb{T}}.\]
where $\K(\mathbb{F}_q)_p$ has trivial $\bT$-action, because $\alpha_1$ is an attaching map in $B\bT$. This has already been proven in \cite[Theorem 3.5]{MR1317117}, so we omit the details. There is a $\bT$-equivariant map of commutative ring spectra 
\[ \THH(\K(\mathbb{F}_q)_p)\rightarrow \K(\mathbb{F}_q)_p \] 
which induces a map of homotopy fixed point spectral sequences and since this map sends $t$ to $t$ and $\alpha_1$ to $\alpha_1$, the differential $d_{2p-2}(t)=t^p\alpha_1$ also occurs in the homotopy fixed point spectral sequence
\[ H^*(\bT,V(1)_*\THH(\K(\mathbb{F}_q)_p))\implies V(1)_*\TC^{-}(\K(\mathbb{F}_q)).\]
Note that we could have also proven this directly by examining $\bT$-equivariant attaching maps in $E\bT$, but for the sake of brevity we give the simpler proof. 
\end{proof}
\begin{cor}\label{cor about d2p-2 diff}
There are differentials $d_{2p-2}(\mu_2)=-t^{p-1}\alpha_1\mu_2$, $d_{2p-2}(\sigma b)=-t^{p-1}\alpha_1 \sigma b$, and $d_{2p-2}(\lambda_1^{\prime})=0$. 
\end{cor}
\begin{proof}
This is immediate from Lemma \ref{d2p-2} and Corollary \ref{perm cycles homotopy fixed}. Note that $\alpha_1\lambda_1^{\prime}=0$ and there are no other classes in the same bidegree as $t^{p-1}\alpha_1\lambda_1^{\prime}$ so $d_{2p-2}(\lambda_1^{\prime})=0$.  
\end{proof}

Now, the classes $\beta_i\in \pi_{(2p^2-2)i-2p}S$ have the property that in the $BP$-Adams spectral sequence for $S$ they are represented by the classes
 \[ \binom{i}{2}v_2^{i-2}k_0 + i v_2^{i-1}b_{1,0} \mod (p,v_1),\]
where $\binom{i}{2}=0$ when $i=1$, which are in $BP$-Adams filtration $2$ (see Definition \ref{beta elements}). We will therefore give a similar argument to the one in the proof of Proposition \ref{perm classes}, except that we will work in the $BP$-Adams spectral sequence in order to use the fact that the classes representing $\beta_k$ are in low $BP$-Adams filtration. To do this we must compute $BP\wedge V(1)_*(T_k(\K(\mathbb{F}_q)_p))$ up to possible differentials of length $4$ or higher.  
 
\subsection{The $BP\wedge V(1)$-Hochschild--May spectral sequence } 
In this section, we begin by computing the input of 
 the $BP\wedge V(1)$-Hochschild--May spectral sequence. 
\begin{lem} \label{inputBPV} There is an isomorphism of $(BP\wedge V(1))_*H\pi_*\K(\bF_q)_p$-algebras
\begin{align}\label{eqBP} \begin{array}{l}(BP\wedge V(1))_*\THH(H\pi_*\K(\bF_q)_p) \cong P(t_1,t_2, \dots )\otimes E(\epsilon_1, \lambda_1)\otimes P(\mu_1)\otimes N ,
\end{array}\end{align}
where $N=\HH_*(S/p_*(H\pi_*\K(\bF_q)_p))$ and the Hurewicz map 
\[ (BP\wedge V(1))_*\THH(H\pi_*\K(\bF_q)_p)\rightarrow (H\bF_p\wedge BP\wedge V(1))_*\THH(H\pi_*\K(\bF_q)_p)\]
sends $t_1$ to $\bar{\xi}_1-\hat{\xi}_1,$ where $\bar{\xi}_1$ is the generator in degree $2p-2$ of $H_*BP$ and $\hat{\xi}_1$ is the generator in degree $2p-2$ of $H_*(V(1)\wedge \THH(H\pi_*\K(\bF_q)_p )$. 
\end{lem}
\begin{proof}
Recall that $V(1)\wedge \THH(H\pi_*( \K(\bF_q)_p))$ is a $V(1)\wedge H\pi_*( \K(\bF_q)_p)$-algebra, and hence an $H\bF_p$ algebra, since $V(1)\wedge H\pi_*( \K(\bF_q)_p)$ is itself an $H\bF_p$-algebra. Thus, there is an equivalence 
\[ BP\wedge  V(1)\wedge \THH(H\pi_*( \K(\bF_q)_p))\simeq BP\wedge H\bF_p\wedge_{H\bF_p} V(1)\wedge \THH(H\pi_*( \K(\bF_q)_p)) \] 
and by the collapse of the K\"unneth spectral sequence and Proposition \ref{computation}, the isomorphism \eqref{eqBP} holds. 

Since $BP\wedge V(1)\wedge \THH(H\pi_*( \K(\bF_q)_p))$ is an $H\bF_p$-module we
can use Lemma \ref{prim}, which states that $(BP\wedge V(1))_*\THH(H\pi_*( \K(\bF_q)_p))$ includes as the comodule primitives inside of 
\[ (H\bF_p\wedge BP \wedge V(1))_*\THH(H\pi_*( \K(\bF_q)_p)). \]
We recall that by the K\"unneth isomorphism and Proposition \ref{Bok1} there is an isomorphism of graded rings 
\begin{align*}
		(H\bF_p\wedge BP \wedge V(1))_*\THH(H\pi_*( \K(\bF_q)_p))\cong Z\otimes  (A//E(0))_*\otimes E(\sigma \bar{\xi}_1)\otimes P(\sigma \bar{\tau}_1)\otimes  N 
\end{align*}
where $N=\HH_*(S/p_*(H\pi_*( \K(\bF_q)_p)))$  and $Z=H_*(BP\wedge V(1))$. We recall from \cite[Lemma 3.4]{K1localsphere} that in degrees less than or equal to $2p-2$ then $N$ has generators 
\[ \{1,x_1,\sigma x_1,x_2\}\]
as a graded $\mathbb{F}_p$-module where $|x_1|=2p-3$, $|\sigma x_1|=2p-2$, and $|x_2|=2p-2$, which will suffice for our purposes. In fact, the elements $1,x_1,\sigma x_1,$ and $x_2$ are also indecomposable algebra generators in $N$ and they are permanent cycles in the algebraic Hochschild--May spectral sequence corresponding to $1,\alpha_1,\sigma \alpha_1,$ and $v_1$ respectively in $\HH_*(S/p_*j)$.

We use the notation
$(\mathcal{A}//E(0))_*\cong P(\hat{\xi}_1,\hat{\xi}_2, \dots ) \otimes E(\hat{\tau}_1,\hat{\tau}_2, \dots )$ 
and 
$H_*(BP)\cong P(\bar{\xi}_1,\bar{\xi}_2,\dots)$ 
to distinguish the two sets of generators. We also recall that  $H_*(V(1))\cong E(\tau_0,\tau_1)$ as a sub $\mathcal{A}_*$-comodule of the dual Steenrod algebra $\mathcal{A}_*$. \footnotemark\footnotetext{As discussed before Proposition \ref{Bok1}, a complete description of $\HH_*(S/p_*(H\pi_*( \K(\bF_q)_p)))$ can be found in \cite[p. 271-272]{K1localsphere}, but we will not need an explicit description here because we will use the algebraic Hochschild--May spectral sequence of \cite[\S 3.1.1]{K1localsphere} to compute the $d^1$ differentials and by the $E^2$-page the description is much simpler.}
The co-action on $\bar{\xi}_i$, $\bar{\tau}_i$, $\hat{\tau}_i$ and $\hat{\xi}_i$ are the same as the co-product in the dual Steenrod algebra, and hence for example $\bar{\xi}_1-\hat{\xi}_1$ is a comodule primitive, since
\[ \psi(\bar{\xi}_1 -\hat{\xi}_1)=1\otimes \bar{\xi}_1 +\bar{\xi}_1\otimes 1 - 1\otimes \hat{\xi}_1  - \bar{\xi}_1\otimes 1 =1\otimes \bar{\xi}_1  - 1\otimes \hat{\xi}_1. \]
 
The co-action on the remaining algebra generators in degrees less than or equal to $2p-2$ is
\[ 
\begin{array}{ll}
\psi(x_1)=1\otimes x_1  \\
\psi(\sigma x_1)=1\otimes \sigma x_1 \\
\psi(x_2)=1\otimes x_2 +\bar{\tau}_0\otimes x_2& 
\end{array} 
\]
and we may observe that there are no other $\mathcal{A}_*$ comodule primitives in degree $2p-2$ other than $\bar{\xi}_1-\hat{\xi}_1$ and $\sigma x_1$ so $t_1$ must map to a linear combination of $\bar{\xi}_1-\hat{\xi}_1$ and $\sigma x_1$. However, the May filtration of $\bar{\xi}_1-\hat{\xi}_1$ is zero and the May filtration of $\sigma x_1$ is positive so, since the Hurewicz map preserves May filtration, we determine that $t_1$ maps to $\bar{\xi}_1-\hat{\xi}_1$.
\end{proof}		

\begin{prop} \label{BPVTHHj}
There is an isomorphism of $(BP_*,BP_*BP)$-comodules
\[ (BP\wedge V(1))_*\THH(\K(\bF_q)_p)\cong P(t_1^p,t_2,...)\otimes E(b)\otimes E(\lambda_1^{\prime},\lambda_2)\otimes P(\mu_2)\otimes \Gamma(\sigma b) \] 
with $BP_*BP$-co-action 
\begin{align*}
		\psi(t_1^p)=& 1\otimes t_1^p+t_1^p\otimes 1, \\
		\psi(\mu_2)=&1\otimes \mu_2, \\ 
		\psi(t_n)=& \sum_{k=0}^nt_k\otimes t_{n-k}^{p^k} \text{ for } n\ge 2,\\ 
		 \psi(\gamma_{n}(\sigma b))=&1\otimes \gamma_{n}(\sigma b), \\
		\psi(b)=&1\otimes b, \\ 
		 \psi(\lambda_2)=&1\otimes \lambda_2+t_1\otimes \lambda_1^{\prime}, \text{ and }\\
		  \psi(\lambda_1^{\prime})=&1\otimes \lambda_1^{\prime}
\end{align*}
where $\lambda_1^{\prime}$ is detected by $\sigma t_1^p$ in the $E^{\infty}$-page, $\lambda_2$ is detected by $\sigma t_2$ in the $E^{\infty}$-page, and $\gamma_{n}(\sigma b)$ is detected by $\gamma_{np}(\sigma \alpha_1)$ in the $E^{\infty}$-page for $k\ge 0$. 
\end{prop}
\begin{proof}
We need to compute differentials in the $BP\wedge V(1)$-Hochschild--May spectral sequence 
\[ E^1_{*,*}=(BP\wedge V(1))_{*,*}\THH(H\pi_*( \K(\bF_q)_p)))\implies BP\wedge V(1)_{*}\THH(\K(\bF_q)_p) \] 
so we examine the map of spectral sequences 
\begin{equation}\label{maposs}
\xymatrix{ (BP\wedge V(1))_{*,*}\THH(H\pi_*( \K(\bF_q)_p)))\ar@{=>}[r]\ar[d]^{h} &  BP\wedge V(1)_{*}\THH(\K(\bF_q)_p) \ar[d] \\
(H\bF_p\wedge BP\wedge V(1))_{*,*}\THH(H\pi_*( \K(\bF_q)_p))\ar@{=>}[r] &  (H\bF_p\wedge BP\wedge V(1))_{*}\THH(\K(\bF_q)_p) . }
\end{equation}
induced by the Hurewicz map $BP\rightarrow H\bF_p\wedge BP.$ Recall from Lemma \ref{inputBPV} that 
\begin{align*}
		(H\bF_p\wedge BP \wedge V(1))_*\THH(H\pi_*( \K(\bF_q)_p))\cong H_*(BP\wedge V(1))\otimes (A//E(0))_*\otimes P(\mu_1)\otimes  N 
\end{align*}
where $N=\HH_*(S/p_*(H\pi_*( \K(\bF_q)_p)))$. We know that in the $H\bF_p\wedge BP\wedge V(1)$-Hochschild--May spectral sequence the classes $\bar{\xi}_i$  for $i\ge 1$ and $\tau_j$ for $j=0,1$ survive to $E^{\infty}$, since the output of the spectral sequence is known to be 
\[ (H\bF_p\wedge BP\wedge V(1))_{*}\THH(\K(\bF_q)_p) \cong  P(t_1,t_2,\dots)\otimes E(\tau_0,\tau_1)\otimes H_*(\THH(\K(\bF_q)_p))\]
by Theorem \ref{HFpj} and the K\"unneth isomorphism. 
This forces the same $d^1$-differentials that occur in the $H\mathbb{F}_p\wedge V(1)$-Hochschild--May spectral sequence and consequently there is an additive isomorphism 
\[  E^2_{*,*}=P(t_1,t_2, \dots )\otimes E(\epsilon_1, \lambda_1,\sigma v_1, \alpha_1)\otimes P(v_1,\mu_1)\otimes \Gamma(\sigma \alpha_1).\]
The map of spectral sequences is therefore again injective on $E^2$-pages. In the $H\mathbb{F}_p\wedge V(1)$-Hochschild--May spectral sequence there are differentials
\begin{align*}
		d^{2p-3}(t_1)=&\alpha_1, \\
		 d^{2p-3}(\lambda_1)=&\sigma \alpha_1,\\
		d^{2p-2}(\hat{\tau}_1)=&v_1, \\ 	
		d^{2p-2}(\mu_1)=&\sigma v_1
\end{align*}
and no further non-trivial differentials besides those differentials generated by these differentials using the Leibniz rule. 
Since the Hurewicz map $h$ is injective and it sends $t_1$ to $\bar{\xi}_1-\hat{\xi}_1$, the differential $d^{2p-3}(t_1)=\alpha_1$ in the top spectral sequence of Equation \eqref{maposs} can be computed using the formula 
\[ d^{2p-3}(h(t_1))=d^{2p-3}(\bar{\xi}_1-\hat{\xi}_1)=\alpha_1=h(\alpha_1).\]
Similarly, $\epsilon_1$ maps to $\bar{\tau}_1-\hat{\tau}_1$ implying $d^{2p-2}(\epsilon_1)=v_1$. 
Hence, in the $BP\wedge V(1)$-Hochschild--May spectral sequence there are differentials 
\[ 
	\begin{array}{llll}
		d^{2p-3}(t_1)=\alpha_1, & d^{2p-3}(\lambda_1)=\sigma \alpha_1, & d^{2p-2}(\epsilon_1)=v_1, & d^{2p-2}(\mu_1)=\sigma v_1 .\\				
	\end{array}
\]
On $E^{2p-1}$-pages the map of spectral sequences induced by the Hurewicz map is again injective. Since $E^{2p-1}\cong E^{\infty}$ in the target spectral sequence, the same is true in the source. This implies that the $BP\wedge V(1)$-Hochschild--May spectral sequence collapses at the $E^{2p-1}$-page. From the $E^{2p-1}$-page on, we write
\begin{align*}
	\gamma_{n}(\sigma b):=&\gamma_{np}(\sigma \alpha_1) \\
	\lambda_1^{\prime}:=&\lambda_1\gamma_{p-1}(\sigma \alpha_1), \\
	\lambda_2:=&(\sigma v_1) \mu_1^{p-1},\\
	\mu_2:=&\mu_1^p, \text{ and }\\
	b := &\alpha_1t_1^{p-1}.
\end{align*}
Write $\mfilt(x)$ for the May filtration of an element $x$. We note that $\mfilt(\gamma_{n}(\sigma b))=(2p-2)np-1$, $\mfilt(\lambda_1^{\prime})=(p-1)(2p-2)-1$, $\mfilt(\lambda_2)=2p-2$ and $\mfilt(b)=2p-3$ for later use. 
By examining the long exact sequence 
\[ BP_*(V(1)\wedge \K(\mathbb{F}_q)_p) \rightarrow BP_*V(1)\wedge \ell \rightarrow BP_*(V(1)\wedge \Sigma^{2p-2}\ell) \] 
we can determine that the co-action on $t_1^p$ and $t_i$ for $i\ge 2$ in $BP_*\K(\mathbb{F}_q)_p$ is the same as the co-action on these elements in $BP_*(V(1)\wedge \ell) \cong P(t_1,t_2,\dots)$. Note that there is no hidden comultiplication on $t_1^p$ since there are no classes in degrees $2p^2-2p-(2p-2)$ or lower and the lowest degree element in $BP_*BP$ is in degree $2p-2$. The class $b$ is the class in lowest positive degree and therefore it is primitive. This produces the co-action on $b,t_1^p,t_i$ for $i\ge 2$ in $BP_*(V(1)\wedge \THH(\K(\mathbb{F}_q)_p))$, by using the splitting of $BP_*BP$-comodules 
\[ BP_*(V(1)\wedge \THH(\K(\mathbb{F}_q)_p)) \cong BP_*(V(1)\wedge \K(\mathbb{F}_q)_p) \oplus BP_*(V(1)\wedge \overline{\THH}(\K(\mathbb{F}_q)_p))\]
induced by the splitting $\THH(\K(\mathbb{F}_q)_p) \simeq \K(\mathbb{F}_q)_p\vee \overline{\THH}(\K(\mathbb{F}_q)_p)$, which we have because $\K(\mathbb{F}_q)_p$ is a commutative ring spectrum. 

The co-action on $\mu_2$ is primitive because $|\mu_2|=2p^2$ and there are no classes in degrees $2p^2-2p+2$ or $2p^2-4p+4$ or lower and the classes in $BP_*BP$ are in degrees congruent to zero mod $2p-2$. Similarly, the co-action on $\lambda_1'$ is primitive because there are no classes in degree $2p^2-2p+1 -(2p-2)$ or lower. 
To determine the co-action on $\lambda_2$, note that there is an isomorphism 
 \[ 
 \begin{array}{rcl} BP_*(V(1)\wedge  \K(\mathbb{F}_q)_p )  &\cong& P(\bar{\xi}_1^p,\bar{\xi}_2, \dots )\otimes E(b) \\
 & \cong  & P(t_1^p,t_2,\dots) \otimes E(b)\end{array} \]
 so $\bar{\xi}_2$ and $t_2$ are two names for the same basis element up to multiplication by a unit. Similarly, $\bar{\xi}_1^p$ and $t_1^p$ are two names for the same basis element up to multiplication by a unit (cf. \cite{MR319197} as described in \cite[p. 511]{MR0458423}). The operation $\sigma$ gives 
 $\lambda_2 = \sigma \bar{\xi}_2\dot{=}\sigma t_2 $ and $\lambda_1'\dot{=}\sigma t_1^p$
 and we can therefore compute the co-action on $\lambda_2$ using the formula
 $\psi(\lambda_2)=(1\otimes \sigma ) (t_2\otimes 1+t_1\otimes t_1^p+1\otimes t_2)$,
 due to \cite{MR2171809}. In other words, in $(BP\wedge V(1))_*\THH(\K\mathbb{F}_q)_p)$, 
 \[ \psi(\lambda_2)\dot{=}1\otimes \lambda_2+t_1\otimes \lambda_1^\prime .\]
This just leaves the classes $\gamma_{n}(\sigma b)$ for $n\ge 1$. Since $\gamma_{p^k}(\sigma b)$ are the indecomposable algebra generators of $\Gamma(\sigma b)$ over a finite field $\mathbb{F}_p$ and products of primitives are primitive, to compute that $\gamma_n(\sigma b)$ is a comodule primitive for each $n\ge 1$ it suffices to show that $\gamma_{p^k}(\sigma b)$ is a comodule primitive for each $k\ge 0$. 
Note that we already checked that in the input of the $BP\wedge V(1)$-Hochschild--May spectral sequence the classes $\gamma_{p^{k+1}}(\sigma \alpha_1)=\gamma_{p^{k}}(\sigma b)$ are primitive. Therefore, it suffices to check that there is not a hidden $BP_*BP$ co-action in the Hochschild--May spectral sequence. 

If the co-action contains terms of the form $x\otimes m$ where $|m|\le |\gamma_{p^k}(\sigma b)|$, then the May filtration of $m$ must be strictly greater than the May filtration of $\gamma_{p^k}(\sigma b)$ (since it would have a to arise from a hidden $BP_*BP$-coaction in the Hochschild--May spectral sequence). Our goal will be to reach a contradiction. 

Suppose that the inequality 
\[\mfilt(m)>\mfilt(\gamma_{p^k}(\sigma b))=(2p-2)p^{k+1}-1\] 
is satisfied. 
Then, since the only classes with positive May filtration are $\gamma_{p^j}(\sigma b)$, $b$, $\lambda_1'$, and $\lambda_2$, the class $m$ must be of the form 
 \[ (\gamma_{\ell}(\sigma b))b^{\epsilon_1} \lambda_1'^{\epsilon_2} \lambda_2^{\epsilon_3} z ,\]
 for some element $z$, where $0\le \ell<p^k$ and $\epsilon_1,\epsilon_2,\epsilon_3\in \{0,1\}$. 
 Write $\mfilt(x)$ for the May filtration of an element, then 
 \begin{align*}
 \mfilt(\gamma_{\ell}(\sigma b))=&(2p-2)p\ell -1\\
 \mfilt(b)=& 2p-3 \\
 \mfilt(\lambda_1')=& (p-1)(2p-2)-1, \text{ and}\\
 \mfilt(\lambda_2)=& 2p-2
\end{align*}
so $\ell,$ $\epsilon_1,$ $\epsilon_2,$  and $\epsilon_3$ must satisfy 
 \begin{align}\label{ineq} 
 (2p-2)p\ell -1+
 \epsilon_1\cdot 2p-3+
 \epsilon_2  \cdot (p-1)(2p-2)-1+ 
 \epsilon_3\cdot  2p-2 \ge (2p-2)p^{k+1}-1. 
 \end{align}
We can rewrite the left-hand side as 
\[  (2p-2)(p\ell +\epsilon_1 + (p-1)\epsilon_2+\epsilon_3)-3\]
and conclude that the inequality 
\[p\ell+\epsilon_1+(p-1)\epsilon_2+\epsilon_3 > p^{k+1}\]
must be satisfied, since $p\ge 5$. Since $p^k>\ell$, we conclude that  
\[ p^{k+1}>p\ell>p^{k+1}-(\epsilon_1+(p-1)\epsilon_2+\epsilon_3)\]
and therefore 
\begin{align}\label{easyinequality}
	\epsilon_1+(p-1)\epsilon_2+\epsilon_3>1.
\end{align}
The degree of $(\gamma_{\ell}(\sigma b)) b^{\epsilon_1} \lambda_1'^{\epsilon_2} \lambda_2^{\epsilon_3}z$ satisfies 
\begin{align*}
| (\gamma_{\ell}(\sigma b)) b^{\epsilon_1} \lambda_1'^{\epsilon_2} \lambda_2^{\epsilon_3}z| >  &(2p^2-2p)(p^{k+1}-(\epsilon_1+(p-1)\epsilon_2+\epsilon_3))\\
&+(2p^2-2p-1)\epsilon_1+(2p^2-2p+1)\epsilon_2+(2p^2-1)\epsilon_3 \\
=& |\gamma_{p^k}(\sigma b) |  -\epsilon_1+((p-2)(2p^2-2p)+1)\epsilon_2+(2p-1)\epsilon_3\\
>& |\gamma_{p^k}(\sigma b) |
\end{align*}
whenever $p\ell  > p^{k+1}-(\epsilon_1+(p-1)\epsilon_2+\epsilon_3)$ and  either $\epsilon_2>0$ or $\epsilon_3>0$.  Since \eqref{easyinequality} cannot be satisfied when $\epsilon_2=\epsilon_3=0$, we have reached a contradiction to the assumption that $|m|<|\gamma_p(\sigma b)|$. Thus, no such $m$ such that $|m|\le |\gamma_{p^k}(\sigma b)|$ and $\mfilt(m) > \mfilt(\gamma_{p^k}(\sigma b))$ exists. This implies that there are no hidden co-actions and $\gamma_{p^k}(\sigma b)$ remains a comodule primitive. 
\end{proof}
\begin{cor} \label{S1BP}
In the homotopy fixed point spectral sequence 
\[ H^*(\bT, (BP\wedge V(1))_*\THH(\K(\mathbb{F}_q)_p))\implies (BP\wedge V(1))^c_*\THH(\K(\mathbb{F}_q)_p)^{h\mathbb{T}} \]
 there are differentials 
\[
	\begin{array}{c} 
		d^2(t_1^p)\dot{=}t\lambda_1' \\
		d^2(t_2)\dot{=}t\lambda_2\\
		d^2(b)\dot{=}t\sigma b 
	\end{array}
\]
and no further $d^2$ differentials besides those generated from these $d^2$ differentials using the Leibniz rule. 
\end{cor}

\begin{proof} 
This follows from Proposition \ref{prop:d2homofixedpoint} and the fact that $\lambda_2\dot{=}\sigma t_2$ and $\lambda_1'\dot{=}\sigma t_1^p$ as discussed in Proposition \ref{BPVTHHj}. 
\end{proof}
 
\begin{rem}\label{rem co-action on t} We also need to know the co-action of $BP_*BP$ on 
\[ BP_*(V(1)\wedge T_{k}(\K(\bF_q)_p)),\] 
which is isomorphic to a sub-quotient of 
\[ P(t_1^p,t_2,...)\otimes E(b)\otimes E(\lambda_1^{\prime},\lambda_2)\otimes P(\mu_2)\otimes \Gamma(\sigma b)\otimes P(t)/t^{k+1} \]
after taking into account the differentials.
Modulo higher filtration, this just amounts to describing the co-action on the class $t$ in the input of the homotopy fixed point spectral sequence 
\begin{align}\label{input of ss} H^*(\bT,BP_*(V(1)\wedge \THH(\K(\bF_q)_p))), \end{align}
since the co-action on a sub-quotient of \eqref{input of ss} is determined by the co-action on \eqref{input of ss}. 

Since we know that $\psi_{H\mathbb{F}_p}(t)=\sum_{i\ge 0} \bar{\xi}_i\otimes t^{p^i}$  
and $t_i\dot{=} \bar{\xi}_i$ in $H_*BP$ (cf. \cite{MR319197} as described in \cite[p. 511]{MR0458423}), 
$BP_*BP$-comodule structure on $H^*(\bT,\mathbb{F}_p)\otimes BP_*/(p,v_1)$ is determined modulo higher filtration by the formula
\[ \psi_{BP}(t)=\sum_{i\ge 0} t_i \otimes t^{p^i} \]
in addition to the usual $BP_*BP$-co-action on $BP_*/(p,v_1)$
\end{rem}
Note that there is a \emph{truncated homotopy fixed point spectral sequence} with signature 
\begin{align}\label{truncated} E^2_{*,*}=P(t)/t^{k+1}\otimes BP_*(V(1)\wedge \THH(\K(\mathbb{F}_q)_p)) \implies BP_*(V(1)\wedge T_{k}(\K(\mathbb{F}_q)_p)).\end{align}
as discussed in Remark \ref{truncated ss}. 
\begin{cor}\label{cor E^4}
In the spectral sequence \eqref{truncated} computing $(BP\wedge V(1))_*T_{k}(\K(\mathbb{F}_q)_p)$, 
there is an isomorphism of $P_{k+1}(t)$-modules and $BP_*BP$ comodules 
\[E^4_{*,*}\cong  M_1\otimes_{P_{k+1}(t)}  M_2\otimes_{P_{k+1}(t)}  M_3 \otimes_{P_{k+1}(t)}  \left (P_{k+1}(t)\otimes P(t_k : k\ge 3)\otimes P(\mu_2) \right )\]
where $M_1=H_*(P(t_1^p)\otimes E(\lambda_1^{\prime})\otimes P_{k+1}(t);d^2)$,  $M_2=H_*(P(t_2)\otimes E(\lambda_2)\otimes P_{k+1}(t);d^2)$, and $M_3=H_*(E(b) \otimes \Gamma(\sigma b) \otimes P_{k+1}(t);d^2)$. Explicitly, we have
\begin{align*} 
M_1=&\mathbb{F}_p\{\lambda_1^{\prime}t_1^{p(j-1)},t^kt_1^{pj}: \nu_p(j)=0, j\ge 1\}\oplus \left ( P_{k+1}(t)\otimes P(t_1^{p^2})\otimes E(\lambda_1^{\prime}t_1^{p^2-p}) \right ),\\
M_2=&\mathbb{F}_p\{\lambda_2t_2^{j-1},t^kt_2^j : \nu_p(j)=0,j\ge 1\}\oplus \left ( P_{k+1}(t)\otimes P(t_2^p)\otimes E(\lambda_2t_2^{p-1}) \right ), \text{ and }\\
M_3=&\mathbb{F}_p\{\gamma_{j}(\sigma b),t^kb\gamma_{j-1}(\sigma b) : \nu_p(j)=0,j\ge 1\}\oplus  \left (P_{k+1}(t)\otimes E(b\gamma_{p-1}(\sigma b))\otimes \Gamma(\gamma_{p}(\sigma b)) \right ). 
\end{align*}
The co-action on an element $x\in E^4_{*,*}$ is determined modulo higher filtration by the co-action of classes in $(BP\wedge V(1))_* \THH(\K(\mathbb{F}_q)_p)$ and the co-action of $t$ from Remark \ref{rem co-action on t} modulo the differential $d^2$ determined in Corollary \ref{S1BP}. The elements in $P(t)/t^k$ and 
\[\mathbb{F}_p\{\lambda_1^{\prime}t_1^{p(j-1)},t^kt_1^{pj_1},\lambda_2t_2^{j-1},t^kt_2^j , \gamma_j(\sigma b),t^kb\gamma_{j-1}(\sigma b) : \nu_p(j)=0, j\ge 1\} \]
are permanent cycles in the truncated homotopy fixed point spectral sequence. 
\end{cor}

\begin{proof} All but the last statement directly follow from Corollary \ref{S1BP}. For the last statement, we note that the homotopy fixed point spectral sequence for the sphere spectrum 
\[ (BP\wedge V(1))_*\otimes P(t) \implies  (BP\wedge V(1))_*^c(\TC^{-}(S))\]
collapses at the $E^2$-page for bidegree reasons. Consequently, the differentials in homotopy fixed point spectral sequence
\[ P(t)\otimes (BP\wedge V(1))_*\THH(R)\implies (BP\wedge V(1))_*^c\TC^{-}(R) \]
are always $t$-linear and the same is true for the truncated homotopy fixed point spectral sequence. This implies that the elements in $P(t)/t^{k+1}$ are permanent cycles and that the elements in 
\[\mathbb{F}_p\{\lambda_1^{\prime}t_1^{p(j-1)},t^kt_1^{pj_1},\lambda_2t_2^{j-1},t^kt_2^j , b\gamma_{j-1}(\sigma b),t^k\gamma_j(\sigma b) : \nu_p(j)=0, j\ge 1\} \]
are permanent cycles.
\end{proof}    

\subsection{Detecting $\beta$ elements in $\TC^{-}$}
First, we define the $\beta$-family. 
\begin{defin}\label{beta elements}
Let $p\ge 5$. By \cite[Example 5.1.20]{rav1}, the elements $\beta_s\in \pi_{(2p^2-2)s-2p}S_p$ are represented by the permanent cycles 
\begin{align}\label{beta} \overline{\beta}_s:=\binom{s}{2}v_2^{s-2}k_0 + sv_2^{s-1}b_{1,0} \mod (p,v_1) \end{align}
in $BP_*\overline{BP}_*^{\otimes_{BP_*}2}$ in the $BP$-Adams spectral sequence where
\begin{align}
\label{k0} k_0:= &2 t_1^p\otimes t_2   -2t_1^p\otimes t_1^{p+1} -t_1^{2p}\otimes t_1\text{ and }\\
\label{b10ANSS} b_{1,0}:= &-\sum_{i=1}^{p-1}\frac{1}{p}\binom{p}{i}t_1^{p-i}\otimes t_1^i . \end{align}
\end{defin}
Recall that $T_{s}(R)$ is defined to be the spectrum $F(S(\bC^{s+1})_+,\THH(R))^{\bT}$. As noted before Proposition 1.4 in \cite{BR05}, $T_{s}(R)$ is a commutative ring spectrum whenever $R$ is a commutative ring spectrum. In particular, $\TC^{-}(R)$ is a commutative ring spectrum. In the following theorem, let $\ell_p$ be a commutative ring spectrum model for the $p$-complete connective Adams summand as constructed in \cite{MR1164148}. Note that there is a map of $E_{\infty}$ ring spectra $L_{K(1)}S\to L_p$  given by inclusion of homotopy fixed points since 
\[L_{K(1)}S\simeq L^{h\mathbb{Z}_p}\] 
where  $\mathbb{Z}_p$-action on $L_p$ is by $E_{\infty}$ ring maps and consequently the inclusion of fixed point $L_{K(1)}S\to L$ is a map of $E_{\infty}$ ring spectra (cf. \cite[Theorem 4.4, Proposition 8.1]{MR2030586}). There is therefore a map of $E_{\infty}$-ring spectra $K(\mathbb{F}_q)_p\to \ell_p$ by applying $\tau_{\ge 0}$ to the map above and using the weak equivalences of $E_{\infty}$-ring spectra $K(\mathbb{F}_q)_p\simeq \tau_{\ge 0}L_{K(1)}S$ and $\tau_{\ge 0}L_p\simeq \ell_p$ where $\tau_{\ge 0}$ is the Postnikov truncation functor (see \cite[p. 260]{K1localsphere} for a discussion of the first weak equivalence  and \cite[Proposition 2.1]{BR08} for the second weak equivalence). 
\begin{prop} \label{AR thm}
The classes $v_2^s$ map to nonzero classes $(t\mu_2)^s$ under the unit map 
\[ V(1)_*S \rightarrow V(1)_*T_{k}(K(\mathbb{F}_q)_p). \]
\end{prop}
\begin{proof} 
By \cite[Prop. 4.8]{MR1947457}, the element $v_2\in V(1)_{2p^2-2}$ maps to $t\mu_2\in V(1)_{2p^2-2}(T_1(\ell_p))$. By \cite[Thm. 6.7]{MR1947457}, the elements $(t\mu_2)^k$ survive the $\mathbb{T}$-homotopy fixed point spectral sequence computing $V(1)_*\TC^{-}(\ell_p)$ and $v_2^s$ maps to $(t\mu_2)^s$ because the map $V(1)_*\to V(1)_*\TC^{-}(\ell_p)$ is a ring map. The map $V(1)_*\TC^{-}(\ell_p)\to V(1)_*T_{k}(\ell_p)$  sends $(t\mu_2)^s$ to $(t\mu_2)^s$ by construction of the $\mathbb{T}$-homotopy fixed point spectral sequence. 

Since we showed $v_2$ maps to $t\mu_2$ under the unit map $V(1)_*S\rightarrow V(1)_*T_1(\K(\mathbb{F}_q)_p)$, the map  $V(1)_*T_1(\K(\mathbb{F}_q)_p)\to V(1)_*T_1(\ell_p)$ sends $t\mu_2$ to $t\mu_2$, and the maps 
\[ V(1)_*S \rightarrow T_{s}(\K(\mathbb{F}_q)_p)\rightarrow V(1)_*T_{s}(\ell_p) \]
are ring maps for $s\ge 1$, the classes $v_2^s$ also map to $(t\mu_2)^s$ under the unit map
\[ V(1)_*S \rightarrow V(1)_*T_{s}(\K(\mathbb{F}_q)_p).\] 
We therefore know that $(t\mu_2)^s$ are permanent cycles in the $BP$-Adams spectral sequence and homotopy fixed point spectral sequences computing $V(1)_*T_{s}\K(\mathbb{F}_q)_p$. 
\end{proof}
We will continue to use notation $d^r$ for differentials in the homotopy fixed point spectral sequence and $d_r$ for differentials in the $BP$-Adams spectral sequence to differentiate the two. 
\begin{thm} \label{beta thm}
Let $p\ge 5$. The elements $\beta_s\in \pi_{(2p^2-2)s-2p}S$ for $s\not \equiv 0\mod p$ map to 
\begin{align}\label{image of betas} \binom{s}{2}(t\mu_2)^{s-2}t^2\lambda_1^{\prime}\lambda_2 + s(t\mu_2)^{s-1} t\sigma b \in V(1)_{(2p^2-2)s-2p}TC^{-}(\K(\bF_q)_p)\end{align}
up to multiplication by a unit in $\mathbb{F}_p$ for certain choices of representatives of the classes $(t\mu_2)^{s-2}t^2\lambda_1^{\prime}\lambda_2$ and $(t\mu_2)^{s-1} t\sigma b$. 
\end{thm}
\begin{proof}
The element $\beta_s$ is represented by the class $\overline{\beta}_{s}$ as in Definition \ref{beta elements} by \cite[Example 5.1.20]{rav1}. We therefore need to check that the classes $\overline{\beta}_s$
map to permanent cycles in the $BP_*BP$-cobar complex for 
$V(1)_*T_{s^{\prime}}(\K(\bF_q)_p).$ From the proof, we will see that we can accomplish this by choosing $s=s^{\prime}$. 

Due to the length of the proof, we break it into several steps:
\begin{enumerate}[labelindent=0pt,labelwidth=\widthof{\ref{last-item}},label=\arabic*.,leftmargin=!]
\item \label{Step 1}
First, we show that $v_2^{s}$ in the $BP_*BP$-cobar complex for $V(1)$ maps to $(t\mu_2)^s$ in the $BP_*BP$-cobar complex for $V(1)_*T_s(\K(\bF_q)_p))$. By Proposition \ref{AR thm}, we know that $v_2^s$ maps to $(t\mu_2)^s$ under the map 
\[ V(1)_*\rightarrow V(1)_*T_s(\K(\mathbb{F}_q)_p).\]
By examining the map of Hochschild--May spectral sequences induced by the unit map $\eta \wedge \id_{V(1)}\colon \thinspace S\wedge V(1)\rightarrow BP\wedge V(1)$ and the subsequent map of homotopy fixed point spectral sequences induced by this same map, we see that $(t\mu_2)^s$ maps to $(t\mu_2)^s$ under the map
\[ V(1)_*T_{s}(\K(\mathbb{F}_q)_p)\rightarrow (BP\wedge V(1))_*T_{s}(\K(\mathbb{F}_q)_p).\]
We also know the map 
\[ \pi_*(\eta \wedge \id_{V(1)})\colon \thinspace \pi_*(S\wedge V(1)) \rightarrow \pi_*(BP\wedge V(1)) \]
sends the class $v_2^s$ to $v_2^s$ since the edge-homomorphism in the $BP$-Adams spectral sequence is a ring homomorphism. We then use the commutative diagram 
\begin{align}\label{square of SS}
\xymatrix{ 
V(1)_*\ar[r] \ar[d] &  V(1)_*T_{s}(\K(\mathbb{F}_q)_p) \ar[d]   \\
BP_*V(1)\ar[r] &  (BP\wedge V(1))_*T_{s}(\K(\mathbb{F}_q)_p)
}
\end{align}
to determine that $v_2^s\in BP_*V(1)$ maps to $(t\mu_2)^s\in(BP\wedge V(1))_*T_{s}(\K(\mathbb{F}_q)_p)$ and also in the map 
of exact couples of the respective $BP$-Adams spectral sequences.\\

\item \label{Step 2}
We now show that the element $\beta_1$ has nontrivial image in $V(1)_*T_{1}(\K(\mathbb{F}_q)_p)$. The element $\beta_1$ is represented by the permanent cycle $b_{1,0}$ in the $E_1$-page of the $BP$-Adams spectral sequence for $V(1)$. This class maps to a class of the same name in $E_1$-page of the $BP$-Adams spectral sequence for $V(1)\wedge T_1(\K(\bF_q)_p)$. The result then follows by Proposition \ref{perm classes} and a result of Zahler \cite{MR319197}, which implies that $b_{1,0}$ maps to $b_0$ via the map from the $BP$-Adams spectral sequence to the Adams spectral sequence (cf. \cite[p. 511]{MR0458423}). \\

\item  \label{Step 3} 
Next we consider the case $s\equiv 1\mod p$. In this case, we know that $\overline{\beta}_{pk+1}=b_{1,0}v_2^{pk}$
represents $\beta_{pk+1}$ in the $BP$-Adams spectral sequence for $V(1)$ for $k\ge 1$. It maps to $b_{1,0}(t\mu_2)^{pk}$
in the $BP$-Adams spectral sequence for $V(1)\wedge T_{pk+1}(\K(\bF_q)_p)$ up to multiplication by a unit. This follows by Step \ref{Step 1} and the fact that the cobar complex for $V(1)\wedge T_{pk+1}(\K(\bF_q)_p)$ is multiplicative.
Since $\overline{\beta}_{pk+1}$ is a permanent cycle in the $BP$-Adams spectral sequence for $V(1)$ for $k\ge 1$ (this follows from \cite[Lemma 5.4]{MR676562}), the class 
$b_{1,0}( t\mu_2)^{pk}$
is an infinite cycle in the $BP$-Adams spectral sequence for $V(1)\wedge T_{pk+1}(\K(\bF_q)_p)$ for $k\ge 1$, but it could still be a boundary. 
It is on the two-line of the $BP$-Adams spectral sequence, so we just need to check that it is not the boundary of a $d^1$ or $d^2$ differential. 
Note that we will prove that, in fact, the element $b_{1,0} (t\mu_2)^{s}$
is never a boundary for any $s\ge 1$. In particular, we also show that $v_2\beta_1$ in the image of the unit map $V(1)\longrightarrow V(1)_*\TC^{-}(K(\mathbb{F}_q)_p)$ since $v_2\beta_1$ is known to be a permanent cycle in $V(1)_*$ by \cite{Tod71}. 
We will break this into two further sub-steps:\\ 

\begin{enumerate}[labelindent=0pt,labelwidth=\widthof{\ref{last-item}},leftmargin=!]
\item \label{Step 3 a} 
If the class  $b_{1,0} (t\mu_2)^{s}$
is the boundary of a $d^1$, then there is a sum of classes
\begin{align}\label{sum possibly hitting b10} \sum_{i} a_i\otimes m_i\in \overline{BP_*BP} \otimes_{BP_*}BP_*(V(1)\wedge T_{s}(\K(\bF_q)_p)) \end{align}
such that 
\[d_1\left (\sum_{i\in I}a_i\otimes m_i \right )= b_{1,0}\otimes ( t\mu_2)^{s}\]
where $s\ge 1$ and $I$ is a finite set. Recall that the co-action on $x\in (BP\wedge V(1))_* T_{s}(\K(\bF_q)_p)$ is of the form 
\[\psi(x)=1\otimes x+\sum_{j\in J} a_j \otimes x_j\] 
where $|x_j|<|x|$ and $J$ is a finite set. 
Observe that the only monomial classes 
$x\in (BP\wedge V(1))_*T_s(\K(\bF_q)_p) $
such that $x=(t\mu_2)^{s}$ or $x_j=(t\mu_2)^{s}$ are classes of the form $(t\mu_2)^{s}y$ for some $y\in (BP\wedge V(1))_*T_s(\K(\bF_q)_p)$ up to multiplication by a unit. The co-action of such a class is
\[\psi(( t\mu_2)^{s}y)=(1\otimes ( t\mu_2)^{s})\psi(y),\] 
so in order to arrange that $x_{j^{\prime}}=(t\mu_2)^{s}$ for some $j^{\prime}$ the coaction $\psi(y)$ must be of the form 
\[ \psi(y)=1\otimes y+ z\otimes 1+ \sum_{k\in K} b_k\otimes y_k\] 
for some finite set $K$.
Since the only classes in 
$ (BP\wedge V(1))_*T_s(\K(\bF_q)_p) $
that have $z\otimes 1$ as a summand in their co-action for some element $z\ne 0,1$ are the classes $t_1^{p}$ and $t_i$ for $i\ge 2$, the class $y$ must be a product $\prod_{i\ge 2}(t_1^p)^{k_1}(t_i)^{k_i}$ where $k_i\ge 0$, $k_i=0$ for all but finitely many $i\ge 2$, and 
$\sum_{i\ge 1}k_i>0$. Since the internal degree of $a_i\otimes ( t\mu_2)^{s}y$ must equal $(2p^2-2)s+2p^2-2p$ in order to be a summand of \eqref{sum possibly hitting b10}, we have
\[|a_i\otimes ( t\mu_2)^{s}y|=(2p^2-2)(s)+|y|+|a_i|=(2p^2-2)s+2p^2-2p\]
so, the degree of $|y|+|a_i|$ must be $2p^2-2p$. However, the class $t_1^p$ is the element of lowest degree of the form $(t_1^p)^{k_1}\prod_{i\ge 2}(t_i)^{k_i}$, where $k_i\ge 0$, $k_i=0$ for all but finitely many $i\ge 1$, and $\sum_{i\ge 1}k_i>0$, and $|t_1^p|=2p^2-2p$. Also, the only classes in degrees less than or equal to $2p^2-2p$ in $BP_*\overline{BP}$ are powers of $v_1$ and $t_1$. Therefore, the only possibility is $a_i=t_1^{p-j}v_1^j$ for some $0\le j <p$ and $y=1$. We know that 
\[ \Delta(t_1^{p-j}v_1^j)=v_1^j\Delta(t_1^{p-j})=v_1^j\cdot (t_1\otimes 1+1\otimes t_1)^{p-j}\]
and we compute 
\[ 
\begin{array}{rcl} 
d_1(t_1^{p-j}v_1^j\otimes ( t\mu_2)^{s}) &=& - \overline{\Delta}(t_1^{p-j}v_1^j)\otimes ( t\mu_2)^{s})\\
&&+ t_1^{p-j}v_1^j\otimes \overline{\psi}(( t\mu_2)^{s})  \\
&=&\overline{\Delta}(t_1^{p-j}v_1^j)\otimes ( t\mu_2)^{s} .\\ 
\end{array}
\]
We observe that there no linear combination of $d_1(t_1^{p-j}v_1^j\otimes ( t\mu_2)^{s})$ that have the element $b_{1,0} (t\mu_2)^s$ as a summand. We conclude that $m_{i^{\prime}}\dot{=}( t\mu_2)^{s}$ for at least one $i^{\prime}$.\\

Now, if $m_{i^{\prime}}\dot{=}( t\mu_2)^{s}$ for only one $i^{\prime}$, then the element $a_{i^{\prime}}$ corresponding to $m_{i^{\prime}}$ must have reduced co-product $b_{1,0}+z$ for some element $z\in BP_*\overline{BP}\otimes_{BP_*}BP_*\overline{BP}$, up to multiplication by a unit; i.e., 
$ \bar{\Delta}(a_{i^{\prime}})\dot{=}b_{1,0}+z .$
The degree of $a_{i^{\prime}}$ must be $2p^2-2p$, so $a_{i^{\prime}}\dot{=}t_1^jv_1^{p-j}$ for $0\le j<p$. 
However, 
\[  
\begin{array}{rcl} 
	\bar{\Delta}(t_1^j v_1^{p-j})&=&  v_1^{p-j}\bar{\Delta}(t_1^j) \\
	&=& v_1^{p-j}(t_1\otimes 1+1\otimes t_1)^j -1\otimes t_1^jv_1^{p-j}-t_1^jv_1^{p-j}\otimes 1
\end{array}   
\]
and this does not equal 
$ b_{1,0}+z$
for any $j$, and any element $z\in BP_*\overline{BP}\otimes_{BP_*}  BP_*\overline{BP} $. 

Suppose that $m_j\dot{=}( t\mu_2)^{s}$ for $j\in J$ where $J$ is a finite set of indices with more than one element.
Then $ \psi(\sum_{j\in J} a_j)  \dot{=} b_{1,0}+z'$
for some element $z'$ in $BP_*\overline{BP}\otimes_{BP_*}  BP_*\overline{BP}$. However, we checked in the proof of Proposition \ref{perm classes} that no class of the form $\sum_{j\in J} a_j$ has co-action $ b_{1,0} +z'$ and the same proof applies here. 

Thus, there is no sum of classes $\sum_{j\in I} a_i\otimes m_i$ such that 
$d_1(\sum_{i\in I}a_i\otimes m_i)\dot{=}b_{1,0} \otimes (t\mu_2)^{s}$
and therefore the class $b_{1,0}(t\mu_2)^{s}$
survives to the $E_2$-page as long as it is a $d^1$ cycle, which is the case when $s=pk$ for $k\ge 1$ and $s=2$ for all $p\ge 5$.\\

\item \label{Step 3 b}
Now suppose there is a class $x$ in bidegree $((2p^2 -2)s+ 2p^2-2p+1,0)$ that is the source of a $d^2$ differential hitting 
$b_{1,0}  (t\mu_2)^{s}.$
Then  
\[x\in BP_{(2p^2 -2)s+ 2p^2-2p+1}(V(1)\wedge T_s(\K(\bF_q)_p))\] 
is a $d_1$-cycle. Since $x$ is in an odd degree, we can classify all the possibilities as a linear combination of elements in the six families of elements, 
 \[ \{  \lambda_1^{\prime}z_1,\lambda_2z_2, t^sbz_3,\lambda_1^{\prime}t_1^{p^2-p}z_4,\lambda_2t_2^{p-1}z_5,b\gamma_{p-1}(\sigma b)z_6 \} \]
 in a sub-quotient of $BP_{(2p^2 -2)s+ 2p^2-2p+1}(V(1)\wedge \THH(\K(\bF_q)_p))\otimes P(t)/t^{s+1}$ by 
 Corollary \ref{cor E^4} where $z_4$, $z_5$, and $z_6$ are non-trivial even dimensional classes and $z_1$, $z_2$, and $z_3$ are non-trivial even dimensional classes that are not $t$-divisible, since $\lambda_1^{\prime}$, $\lambda_2$, and $t^sb$ are simple $t$-torsion elements. 
We know that each of the elements in the set 
\[\{ \lambda_2 , \lambda_1^{\prime}, t^sb\}\] 
are permanent cycles in the truncated homotopy fixed point spectral sequence computing $(BP\wedge V(1))_*T_s(\K(\bF_q)_p)$  by Corollary \ref{cor E^4}. 
We know that each of the elements in the set
\[\{ b\gamma_{p-1}(\sigma b), \lambda_2t_2^{p-1},\lambda_1t_1^{p^2-p}\} \]
are permanent cycles by mapping to a homological homotopy fixed point spectral sequence via the Hurewicz map and applying \cite[Theorem 5.1 (b)]{BR05}. 
Note that the classes $t^sb$, $\lambda_1^{\prime}t_1^{p^2-p}$, $b\gamma_{p-1}(\sigma b)$, and $\lambda_2t_2^{p-1}$ are indecomposable in $BP_{(2p^2 -2)s+ 2p^2-2p+1}(V(1)\wedge T_s(\K(\bF_q)_p))$ by Corollary \ref{cor E^4}.\\
 
We explicitly compute 
 $d_1(\lambda_2)=\bar{\xi}_1\otimes \lambda_1'$. Therefore, by the Leibniz rule, 
$d_1(\lambda_2z_2)=(\bar{\xi}_1\otimes \lambda_1')z_2 + \lambda_2d_1( z_2)\ne 0.$
 Therefore, the classes  of the form $\lambda_2z_2$ do not survive to the $E_2$ page of the $BP$-Adams spectral sequence and cannot be the source of a $d_2$ differential hitting 
 $b_{1,0} (t\mu_2)^{s}.$ \\

Next, we check elements of the form $t^sbz_3$ or $\lambda_1'z_1$ where $z_3$ is not divisible by $t\mu_2$ or $t\sigma b$ and $z_1$ is not divisible by $\lambda_1^{\prime}$. Note that the Leibniz rule implies
\[ d_2(t^sbz_3)=d_2(t^sb)z_3+t^sbd_2(z_3)\]
and similarly, 
\[ d_2(\lambda_1^{\prime}z_2)=d_2(\lambda_1^{\prime})z_2-\lambda_1^{\prime}d_2(z_2)\]
so we need to check if 
\[\alpha(d_2(tb)z_3+tbd_2(z_3))+\beta(d_2(\lambda_1')z_2-\lambda_1'd_2(z_2))=b_{1,0}(t\mu_2)^s\]
for some $\alpha,\beta\in \mathbb{F}_p$. Suppose $d_2(\lambda_1^{\prime}z_2)\dot{=}b_{1,0} (t\mu_2)^s+z$ for some possibly trivial class $z$. Then by the Leibniz rule 
$d_2(\lambda_1^{\prime}z_2)=d_2(\lambda_1^{\prime})z_2-\lambda_1^{\prime}d_2(z_2)$
and so $d_2(\lambda_1^{\prime})$ would need to divide the product of $b_{1,0}$ and $(t\mu_2)^s$. Since we already proved that $b_{1,0}$ and $(t\mu_2)^s$ are permanent cycles, this would be a contradiction. It therefore suffices to show that $d_2(t^sbz_3)\ne b_{1,0}(t\mu_2)^s+z^{\prime}$ for some class $z^{\prime}$. Since $t\mu_2$ does not divide $z_3$ it suffices to show that $d_2(z_3)t^sb \ne  b_{1,0} (t\mu_2)^s+z^{\prime}$ for some class $z^{\prime}$ by the Leibniz rule. This follows because $t^sb$ does not divide $b_{1,0}(t\mu_2)^s$.\\ 

Finally, we need to check the classes $\lambda_1^{\prime}t_1^{p^2-p}z_4$, $\lambda_2t_2^{p-1}z_5$, and  $b\gamma_{p-1}(\sigma b)z_6$. In order to be a source of a $d_2$, these classes must survive to the $E_2$-page, so we just need to check this in the case when the classes of the form $\lambda_1^{\prime}t_1^{p^2-p}z_4$, $\lambda_2t_2^{p-1}z_5$, and  $b\gamma_{p-1}(\sigma b)z_6$ are $d_1$ cycles. Since the $d_1$-cycles on the $0$-line are exactly that $BP_*BP$-comodule primitives, we just need to 
prove that the elements 
\begin{align}\label{some elements} \{ \lambda_1^{\prime}t_1^{p^2-p}, b\gamma_{p-1}, \lambda_1^{\prime}t_1^{p^2-p}\}\end{align}
are not $BP_*BP$-comodule primitives. From Corollary \ref{cor E^4}, we know that they are not $BP_*BP$-comodule primitives on the $E^4$-page of the truncated homotopy fixed point spectral sequence, which is the $E^{\infty}$-page for $BP_*(V(1)\wedge T_3(\K\mathbb{F}_q)_p))$. We also claim that hidden $BP_*BP$-comodule extensions in the truncated homotopy fixed point spectral sequence computing $BP_*V(1)\wedge T_s(\K(\mathbb{F}_q)_p)$ cannot produce additional comodule primitives. 
To see this, note that elements in higher skeletal filtration are always divisible by higher powers of $t$ and therefore summing with such classes cannot reduce the number of terms in the $BP_*BP$-coaction of an element in the $E^{\infty}$-page in lower skeletal filtration. Therefore, the classes $\lambda_1^{\prime}t_1^{p^2-p}z_4$,  $\lambda_2t_2^{p-1}z_5$, and $b\gamma_{p-1}(\sigma b)z_6$ do not survive to the $E_1$-page of the $BP$-Adams spectral sequence for $BP_*V(1)\wedge T_3(\K(\mathbb{F}_q)_p)$. The elements in \eqref{some elements} in the $BP$-Adams spectral sequence for $V(1)\wedge T_s(\K(\mathbb{F}_q))$ map to classes of the same name in the $BP$-Adams spectral sequence via the map 
 $E_1$-pages of $BP$-Adams spectral sequences 
\[ BP_*\overline{BP}_*^{\otimes_{BP_*} \bullet}  \otimes BP_*V(1)\wedge T_s(\K(\mathbb{F}_q)) \longrightarrow BP_*\overline{BP}_*^{\otimes_{BP_*} \bullet}\otimes BP_*V(1)\wedge T_3(\K(\mathbb{F}_q)) \]
induced by the canonical map $V(1)\wedge T_s(\K(\mathbb{F}_q))\to V(1)\wedge T_3(\K(\mathbb{F}_q))$ and therefore they cannot be comodule primitives in the source, because this would lead to a contradiction. 
 
In fact, our proof shows that $b_{1,0} (t\mu_2)^s+z^{\prime\prime}$ is not hit by a $d_2$ differential for any $s\ge 0$ and any class $z^{\prime\prime}$, so the exact same argument implies that, as long as $b_{1,0}(t\mu_2)^s$ and $\overline{\beta}_s$ are $d_1$-cycles and they are not $d_1$ boundaries, then they are not $d_2$ boundaries for any $s\ge 1$. Note that by \cite[Theorem 5.2]{Tod71}, the element $v_2\beta_1$ is non-trivial in $V(1)_{(2p^2-2)+(2p^2-2p-2)}$. Our argument also shows that $v_2\beta_1$ has non-trivial image in $V(1)_*\TC^{-}(K(\mathbb{F}_q)_p)$. We also use this to reduce the work we need to do to show that $\overline{\beta}_s$ is a permanent cycle in Step \ref{Step 5}
\end{enumerate}

\item \label{Step 4}
We now discuss how to detect $\beta_{2}$. The element $\beta_2$ is represented by the class $k_0+2b_{1,0}v_2 \mod (p,v_1)$ in the input of the $BP$-Adams spectral sequence for $S$. It is also a permanent cycle in the $BP$-Adams spectral sequence for $V(1)$ as a consequence of \cite[Lemma 5.4]{MR676562}. It maps to the class $k_0+2b_{1,0}t\mu_2$ in 
\[  BP_*\overline{BP}_*\otimes_{BP_*} BP_*\overline{BP}_* \otimes_{BP_*} BP_*V(1)\wedge T_2(\K(\mathbb{F}_q)_p) \]
under the map of $E_1$-pages of $BP$-Adams spectral sequences induced by the map $V(1)\rightarrow V(1)\wedge T_2(\K(\mathbb{F}_q)_p)$, by Step \ref{Step 1} and  the multiplicativity of $E_1$-page of the $BP$-Adams spectral sequence. 

Recall that 
$BP_*(V(1)\wedge T_2(\K(\mathbb{F}_q)_p)) $
is isomorphic to a sub-quotient of 
\[ H_*\left(P(t_1^p,t_2, \dots)\otimes E(b)\otimes E(\sigma t_1^p,\sigma t_2)\otimes P(\mu_2)\otimes \Gamma(\sigma b)\otimes P_3(t) ;  d^2(x)=t\sigma x\right)\] 
after taking into account $d^4$ differentials. 
We can therefore check every element in degree $4p^2-2p-2$ in 
\[ BP_*\overline{BP}_*\otimes_{BP_*} BP_*(V(1)\wedge T_2(\K(\mathbb{F}_q)_p)) \]
to see if it has $k_0+2b_{1,0}v_2$ as a boundary up to multiplication by a unit. Note that the source of a $d_1$ differential hitting $k_0+2b_{1,0}v_2$ must be in bidegree $(1,4p^2-2p-2)$ or in other words stem $4p^2-2p-3$ and $BP$-Adams filtration $1$ because the $d_1$ differential decreases stem by one and increases $BP$-Adams filtration by one. Consequently, the source of a possible differential of length $1$ hitting $k_0+2b_{1,0}v_2$ must be of the form $x\otimes m$ where $0\le |x|\le 4p^2-2p-2$ and $0\le |m|\le 4p^2-2p-2$. The elements $x\in BP_*\overline{BP}_*$ satisfying $0\le |x|\le 4p^2-2p-2$ are linear combinations of $t_1^k$ and $t_2$ for $1\le k  < 2p$. We conclude that $4p^2-2p-2 - |m| \equiv 0 \mod 2p-2$. By Corollary \ref{cor E^4}, the potential elements $m\in BP_*V(1)\wedge T_2(\K(\mathbb{F}_q)_p)$ satisfying  $0\le |m|\le 4p^2-2p-2$ are 
\begin{align*}
\{ 1, \lambda_1^{\prime},\lambda_2,\sigma b, \gamma_2(\sigma b), \mu, t\mu, t^2\mu, t^2t_1^p,t^2t_1^{2p},t^2t_2,t^2b,t^2t_1^pb,t^2bt_2\} 
\end{align*}
so after removing all elements that don't satisfy $4p^2-2p-2 - |m| \equiv 0 \mod 2p-2$, we are left with the elements
\[ \{ \sigma b, \gamma_2(\sigma b),t\mu_2 \}.\]

When $m=\sigma b$, the only possibilities for $x$ are linear combinations of $t_1^{p+1}$ and $t_2$, so we need to compute $d_1(t_1^{p+1}\otimes \sigma b)$ and $d_1(t_2\otimes \sigma b)$. When $m=\gamma_{2}(\sigma b)$ the only option is $x=t_1$, so we need to compute $d_1(t_1\otimes \gamma_2(\sigma b))$. When $m=t\mu_2$, $x$ must be $t_1^p$ so we need to compute $d_1(t_1^p\otimes t\mu_2)$. We compute 
\begin{align*}
d_1(t_1^{p+1}\otimes \sigma b) =& -\sum_{i=1}^{p} \binom{p+1}{i} t_1^i\otimes t_1^{p+1-i}\otimes \sigma b\\
d_1(t_2\otimes \sigma b) =&-t_1^{p}\otimes t_1\sigma b \\
d_1(t_1\otimes \gamma_2(\sigma b))=&0\\
d_1(t_1^p\otimes t\mu_2)=&-\sum_{i=1}^{p-1}\binom{p}{i}t_1^i\otimes t_1^{p-i}\otimes t\mu=0
\end{align*} 
We then note that  the equation 
\[ a_1d_1(t_1^{p+1}\otimes \sigma b)+a_2d_1(t_2\otimes \sigma b) +a_3d_1(t_1^p\otimes t\mu_2) = k_0+2b_{1,0}v_2\]
is not satisfied for any $a_1,a_2,a_3\in \mathbb{F}_p$. We conclude that $k_0+2b_{1,0}v_2$ survives to the $E_2$-page of the $BP$-Adams spectral sequence.

We next need to check if it is the boundary of a $d_2$. 
By Corollary \ref{cor E^4}, we know that $BP_{4p^2-4p-5}(V(1)\wedge T_2(\K(\mathbb{F}_q)_p))=\mathbb{F}_p\{t^2t_1b\}$. However, we know $\psi(t_1)=t_1\otimes 1+1\otimes t_1$, so $d_1(t_1)=t_1\otimes 1$. Since $t^2b$ is a comodule primitive, we compute $d_1(t_1t^2b)=t_1\otimes t^2b$ and consequently, $t_1t^2b$ does not survive to the $E_2$-page of the $BP$-Adams spectral sequence. Therefore, the class $t_1t^2b$ cannot be the source of a $d_2$ differential hitting $k_0+2b_{1,0}v_2$. 
Consequently, the class $k_0+2b_{1,0}v_2$ survives to the $E^3$-page of the $BP$-Adams spectral sequence. Since it is in $BP$-Adams filtration $2$ and it is an infinite cycle, it also survives to the $E_{\infty}$-page of the $BP$-Adams spectral sequence. 

We now identify an element in $V(1)_{4p^2-2p-4}T_2(\K(\mathbb{F}_q)_p)$ which projects onto $k_0+2b_{1,0}t\mu_2$ in the $E_{\infty}$-page of the $BP$-Adams spectral sequence. 
We claim that the only classes in $V(1)_{4p^2-2p-4}T_2(\K(\mathbb{F}_q)_p)$ are $t^2\lambda_1'\lambda_2$ and $t^2\mu_2\sigma b$. To see this, note that the there is a truncated homotopy fixed point spectral sequence 
\[ P_2(t)\otimes V(1)_*\THH(\K(\mathbb{F}_q)_p)\implies V(1)_{4p^2-2p-4}T_2(\K(\mathbb{F}_q)_p)\]
where the input is
\[ \mathbb{F}_p\{ 1, \alpha_1,\lambda_1^{\prime}, \alpha_1\lambda_2,\lambda_1^{\prime}\lambda_2,\alpha_1\lambda_1^{\prime}\lambda_2\}\otimes P(\mu_2)\otimes \Gamma(\sigma b)\otimes P_2(t).\]
Concentrating on elements in degree $4p^2-2p-4$ allows us to remove all elements in odd degree and degree greater than $4p^2-2p-4$ from consideration, leaving us with considering elements in 
\[ \mathbb{F}_p\{ 1, \alpha_1\lambda_2,\lambda_1^{\prime}\lambda_2\}\otimes E(\mu_2)\otimes P_2(\sigma b)\otimes P_2(t).\]
Direct computation shows that the only possible elements in the correct degree in the $E^2$-page of the truncated homotopy fixed point spectral sequence are $t^2\lambda_1^{\prime}\lambda_2$ and $t^2\mu \sigma b$. 
We already showed that $v_2\beta_1$ has nontrivial image in $V(1)_{4p^2-2p-4}\TC^{-}(K(\mathbb{F}_q)_p)$ in the previous step. We also just showed that $\beta_2$ has nontrivial image. Therefore, both $t^2\lambda_1^{\prime}\lambda_2$ and $t^2\mu_2\sigma b$ must be permanent cycles or else we reach a contradiction. Since we know $v_2$ maps to $t\mu_2$ and $\beta_1$ maps to $t\sigma b$, multiplicativity implies that $v_2\beta_1$ is detected by $t^2\mu_2\sigma b$. 
Therefore, we conclude that the Hurewicz image of $\beta_2$ is detected by a linear combination 
$ct^2\lambda_1^{\prime}\lambda_2+c^{\prime}t^2\mu_2\sigma b\in V(1)_{4p^2-2p-4}T_2(\K(\mathbb{F}_q)_p)$ where $c\ne 0$. 
We may arrange that $\beta_2$ maps to $t^2\lambda_2\lambda_1^{\prime}+2\mu_2\sigma b$ by judiciously choosing our choice of representatives for the classes $t^2\mu_2\sigma b$ and $t^2\lambda_2\lambda_1^{\prime}$ and our basis for the vector space $\mathbb{F}_p\{t^2\mu_2\sigma b,t^2\lambda_2\lambda_1^{\prime}\}$. In fact, we can conclude that $k_0$ representing $\beta_2-2v_2\beta_1$ and $v_2b_{1,0}$ representing $v_2\beta$ also have non-trivial Hurewicz image. 

\item \label{Step 5}
Finally, discuss how to detect $\beta_s$ where $s\not \equiv 0,1 \mod p$. In this case, $\beta_s$ is represented by $\overline{\beta}_s$, which maps to 
\[(t\mu_2)^{s-1}(\binom{s}{2} k_0+kb_{1,0}(t\mu_2))\] 
by Step \ref{Step 1}
We just need to check the $d_1$ differentials on classes of the form $(t\mu_2)^{s-1}w$ where $w$ is an element in $BP_*\overline{BP}_*\otimes_{BP_*}BP_*\overline{BP}_*$. 
If the class  
\[2t_1^p\otimes t_2\otimes ( t\mu_2)^{s-1} -2t_1^p\otimes t_1^{1+p}\otimes ( t\mu_2)^{s-1} -t_1^{2p}\otimes  t_1\otimes ( t\mu_2)^{s-1} \]
is the boundary of a $d_1$, then there is a sum of classes
\[\sum_{i\in I} a_i\otimes m_i\in \overline{BP_*BP }\otimes_{BP_*}BP_*V(1)\wedge T_{k}(\K(\bF_q)_p)\]
for some finite set $I$ such that 
\[d_1(\sum_{i\in I}a_i\otimes m_i)\dot{=} 2t_1^p\otimes t_2\otimes (t\mu_2)^{s-1} -2t_1^p\otimes t_1^{1+p}\otimes  (t\mu_2)^{s-1}-t_1^{2p}\otimes  t_1\otimes  (t\mu_2)^{s-1}.\]

\begin{enumerate}[labelindent=0pt,labelwidth=\widthof{\ref{last-item}},leftmargin=!]
\item 
Recall that the co-action on $x$ is of the form 
\[\psi(m)=1\otimes x+\sum_{j\in J} a_j \otimes x_j\]
for some finite set $J$ where $|x_j|<|x|$. Again, observe that the only elements in 
$(BP\wedge V(1))_*T_k(\K(\bF_q)_p) $
whose co-action contains $( t\mu_2)^{s-1}$ as either $x$ or $x_j$ for some $j$ are classes of the form $(t\mu_2)^{s-1}y$ for some $y\in (BP\wedge V(1))_*T_s(\K(\bF_q)_p)$ up to multiplication by a unit. The co-action of such a class is
$\psi(( t\mu_2)^{s-1}y)=(1\otimes ( t\mu_2)^{s-1})\psi(y),$
and $\psi(y)$ must be of the form 
$\psi(y)=1\otimes y+ z\otimes 1+ \sum b_i\otimes y_i$
since $\psi(( t\mu_2)^{s-1}y)$ must have $1\otimes (t\mu_2)^{s-1}$ as a term, up to multiplication by a unit. Since the only classes in 
$(BP\wedge V(1))_*T_s(\K(\bF_q)_p) $
that have a term $z\otimes 1$ in their co-action are the classes $t_1^{p}$, $t_i$ for $i\ge 2$ the class $y$ must be a product of these. Since $|(u\cdot t\mu_2)^{s-1}y|=(2p^2-2)(s-1)+|y|$ and the degree must equal $4p^2-2p+s(2p^2-2)$, 
the degree of $y$ must be $4p^2-2p$. However, the class $t_1^p$ is the element of lowest degree of the form 
\[t_1^{pi_1}\prod_{j\ge 2}t_j^{i_j}\] 
where $i_j\ge 0$, $i_j=0$ for all but finitely many $j\ge 1$, and $\sum_{j\ge 2} i_j>0$. Also, $|t_1^p|=2p^2-2p$ and the next lowest degree element is $t_2$ with $|t_2|=2p^2-2$ so no product of classes in this set can be in degree $4p^2-2p$. Thus, $m_{i^{\prime}}\dot{=} (t\mu_2)^{s-1}$ for at least one $i^{\prime}\in I$. 

Now, if $m_{i^{\prime}}\dot{=}( t\mu_2)^{s-1}$ for only one $i^{\prime}$, then the element $a_{i^{\prime}}$ corresponding to $m_{i^{\prime}}$ must have reduced co-product 
$2t_1^p\otimes t_2 -2t_1^p\otimes t_1^{1+p}\otimes 1-t_1^{2p}\otimes  t_1$; i.e 
\[ \overline{\Delta}(a_{i^{\prime}})\dot{=} 2t_1^p\otimes t_2 -2t_1^p\otimes t_1^{1+p}-t_1^{2p}\otimes  t_1.\] 
The degree of $a_{i^{\prime}}$ must be $4p^2-2p$, so $a_{i^{\prime}}\dot{=} t_1^jv_1^{p-j}v_2^{\epsilon_1}t_2^{\epsilon_2}$.  
 
However, 
\begin{align*}
	\Delta(t_1^j v_1^{p-j}v_2^{\epsilon_1}t_2^{\epsilon_2})=& v_1^{p-j}v_2^{\epsilon_1} (t_1\otimes 1 +1\otimes t_1)^j(t_2\otimes 1+ 1\otimes t_2+t_1^p\otimes t_1 )^{\epsilon_2}
\end{align*}
and so $\overline{\Delta}(t_1^j v_1^{p-j}v_2^{\epsilon_1}t_2^{\epsilon_2})$ does not equal 
\[ t_1^p\otimes t_2 -2t_1^p\otimes t_1^{1+p}-t_1^{2p}\otimes  t_1\]
for any $0\le j\le p$.\\ 

Suppose that $m_k=( t\mu_2)^{s-1}$ for $k\in K\subset I$ where $K$ contains more than one element. Then 
\[ \psi(\sum_{k\in K} a_k)  \dot{=} 2t_1^p\otimes t_2 -2t_1^p\otimes t_1^{1+p}-t_1^{2p}\otimes  t_1\]
for some possibly trivial element $z^{\prime}\in BP_*\overline{BP}_*\otimes_{BP_*}BP_*\overline{BP}_*$. However, we checked in Step \ref{Step 4} that no class of the form $\sum_{k\in K} a_k$ has co-action 
\[  (2t_1^p\otimes t_2 -2t_1^p\otimes t_1^{1+p}-t_1^{2p}\otimes  t_1) \]
and the same proof applies here. 

Thus, there is no sum of classes $\sum_{i\in I} a_i\otimes m_i$ such that 
\[d_1(\sum_{i\in I} a_i\otimes m_i)= (2t_1^p\otimes t_2 -2t_1^p\otimes t_1^{1+p}-t_1^{2p}\otimes t_1)\otimes (t\mu_2)^{s-1} \]
and therefore the class 
$ (2t_1^p\otimes t_2 -2t_1^p\otimes t_1^{1+p}-t_1^{2p}\otimes t_1)\otimes (t\mu_2)^{s-1} $
survives to the $E_2$-page.\\

\item The argument that there are no $d_2$ differentials that hit 
\[\binom{s}{2}(t\mu_2)^{s-1}k_0+sb_{1,0}(t\mu_2)^{s}\] 
is the same as that of Step \ref{Step 3 b}.

Therefore, $\binom{s}{2}(t\mu_2)^{s-1}k_0+sb_{1,0}(t\mu_2)^{s}$  is a permanent cycle in the $BP$-Adams spectral sequence for $V(1)\wedge T_s(\K(\mathbb{F}_q)_p)$. 

We now just need to show that 
\[ \binom{s}{2}(t\mu_2)^{s-1}(t\lambda_1'\lambda_2)+s(t\sigma b)(t\mu_2)^{s} \in V(1)_{(2p^2-2)s-2p}T_s(\K(\mathbb{F}_q)_p)\]
projects onto the permanent cycle $\binom{s}{2}(t\mu_2)^{s-1}k_0+sb_{1,0}(t\mu_2)^{s}$ in the $E_{\infty}$-page of the $BP$-Adams spectral sequence. Since $\beta_2-2v_2\beta_1$, $\beta_1$, and $v_2$ in $V(1)_*$ map to $t^2\lambda_1^{\prime}\lambda_2$, $t\sigma b$, and $t\mu_2$ by the previous steps and these project onto $k_0$, $\beta_1$, and $v_2$ respectively in the $E_{\infty}$-page of the $BP$-Adams spectral sequence, we can conclude that  $\binom{s}{2}(t\mu_2)^{s-1}(t\lambda_1'\lambda_2)+s(t\sigma b)(t\mu_2)^{s}$ projects onto $\binom{s}{2}(t\mu_2)^{s-1}k_0+sb_{1,0}(t\mu_2)^{s}$ in the $E_{\infty}$-page of the $BP$-Adams spectral sequence because the unit map $V(1)\longrightarrow V(1)\wedge T_sK(\mathbb{F}_q)_p)$ induces a multiplicative map of $BP$-Adams spectral sequences.
\end{enumerate}
\end{enumerate}
\end{proof}

\subsection{Detecting $\beta$ elements in algebraic K-theory}
The goal of this section is to prove that the $\beta$-family mod $(p,v_1)$ is detected in the iterated algebraic K-theory of finite fields. We prove this as a Corollary to Theorem \ref{beta thm}. The proof relies on the fact that the trace map $\K(R)\rightarrow \TC^{-}(R)$ is a map of commutative ring spectra when $R$ is a commutative ring spectrum. The proof that the cyclotomic trace map $\K(R)\to \TC(R)$ is a map of commutative ring spectra when $R$ is a commutative ring spectrum is due to Hesselholt--Geisser \cite{MR1743237} for Eilenberg-MacLane spectra and Blumberg--Gepner--Tabuada \cite{BGT14} for commutative ring spectra. The advantage of the approach of Blumberg--Gepner--Tabuada \cite{BGT14} is that they prove that the multiplicative cyclotomic trace map is also unique. This work builds on their proof that algebraic K-theory is the universal additive functor \cite{MR3070515} (for a different take on universal properties of algebraic K-theory see Barwick \cite{MR3465850}). 

Let $R$ be a connective commutative ring spectrum. We write $\TP(R):=\THH(R)^{t\bT}$. Recall from \cite[Prop. II.1.9]{NS18} and \cite[Rem. II.4.3]{NS18}, that $\TC$ is the homotopy equalizer of the maps 
\begin{align*}
	\text{can} \colon \thinspace \TC^{-}(R)\to \TP(R)^{\wedge}=\prod_{p} \TP(R)_p\\
	\varphi_p^{h\mathbb{T}} \colon \thinspace \TC^{-}(R)\to \TP(R)^{\wedge}=\prod_{p} \TP(R)_p
\end{align*}
where $\text{can} \colon \thinspace \TC^{-}(R)\to \TP(R)^{\wedge}\to \TP(R)$ is the canonical map from $\TC^{-}(R)$ to the cofiber of the norm map $\Sigma\THH(R)_{h\mathbb{T}}\to \TC^{-}(R)$, and 
\[\varphi_p\colon \thinspace \THH(R)\to \THH(R)^{tC_p}\] 
is the Tate valued Frobenius map. This description relies on the identification 
\[ (\THH(R)^{tC_p})^{h\mathbb{T}}\simeq \TP(R)_p\]
from \cite[Lemma II.4.2.]{NS18}, which uses the assumption that $R$ is bounded below. Since the maps $\text{can}$ and $\varphi_p^{h\mathbb{T}}$ are each commutative ring spectrum maps and the forgetful functor $\Comm(\Sp)\to \Sp$ is a right adjoint, we can take the equalizer of $\text{can}$ and $\varphi_p^{h\mathbb{T}}$ in $\Comm(\Sp)$. Therefore, the map $\TC(R)\to \TC^{-}(R)$ is a map of commutative ring spectra. The following is then an easy consequence of \cite[Thm. 7.4]{BGT14}. The author would like to thank the anonymous referee for pointing out this simpler argument.

\begin{lem}
Suppose $R$ is a connective commutative ring spectrum, then the trace map 
\[ \K(R)\to \TC(R)\to \TC^{-}(R)\]
is a map of commutative ring spectra.
\end{lem}

\begin{cor}\label{betafamilyiniteratedKtheory}
Let $p\ge 5$ be a prime. Then each element of the $p$-primary $\beta$-family $\{\beta_k : k\not \equiv 0 \mod p\}$ maps to a nonzero element in $\pi_*\K(\K(\bF_q))$ under the unit map. If $q=\ell^m$ for some prime $\ell$, then each element of the $p$-primary $\beta$-family $\{\beta_k : k\not \equiv 0 \mod p\}$ maps to a nonzero element in $\pi_*\K(\K(\bF_\ell))$ under the unit map.
\end{cor}
\begin{proof}
The elements $\beta_k\in \pi_{(2p^2-2)k-2p}(S)$ for $k\not\equiv 0\mod p$ map non-trivially to 
\[V(1)_{(2p^2-2)k-2p}\TC^{-}(\K(\bF_q)_p)\] under the unit map by Theorem \ref{beta thm}. 
Since the cyclotomic trace is multiplicative, the unit map
\[ \pi_*S\rightarrow V(1)_*\TC^{-}(\K(\bF_q)_p)\]
factors through $\pi_*\K(\K(\bF_q))$; i.e, there is a commutative diagram 
\[ 
	\xymatrix{ 
	  	\pi_*S \ar[d]^{\pi_*\eta_{V(1)}\wedge S} \ar[rr]^{\pi_*\eta_{\K(\K(\mathbb{F}_\ell))}}  &&\pi_*\K(\K(\mathbb{F}_\ell))  \ar[r]^{h} &   \pi_*\K(\K(\mathbb{F}_q)) \ar[d]^{g} \\
		V(1)_*S \ar[drrr]_(.4){V(1)_*\eta_{\TC^{-}(\K(\bF_q)_p)}}  \ar[rrr]^(.4){V(1)_*\eta_{\K(\K(\bF_q)_p)}}  &&& V(1)_*\K(\K(\bF_q)_p) \ar[d]^{V(1)_* tr}  \\
						& &&  V(1)_*\TC^{-}(\K(\bF_q)_p) 
	}
\]
where $g$ is the composite of the unit map $\pi_*\eta_{V(1)}\wedge \K(\K(\mathbb{F}_q) \colon \thinspace \pi_*\K(\K(\mathbb{F}_q))\longrightarrow  V(1)_*\K(\K(\mathbb{F}_q))$ with the map $V(1)_*K(f_p)$ where $f_p$ is the $p$-completion map $f_p\colon\thinspace \K(\mathbb{F}_q)\longrightarrow K(\mathbb{F}_q)_p$.  The map $h$ is induced by the canonical map of fields $\mathbb{F}_{\ell}\longrightarrow \mathbb{F}_q$. 
\end{proof}
\begin{cor}\label{betafamilyiniteratedKtheory2}
Let $R$ be any $E_2$ ring spectrum equipped with a map of $E_2$ ring spectra to $H\mathbb{F}_q$ and let $p\ge 5$. 
Then the $p$-primary $\beta$-family $\{\beta_s : s\not \equiv 0 \mod p\}$ is detected in $\K(\K(R))$. In particular, the $p$-primary $\beta$-family $\{\beta_s : s\not \equiv 0 \mod p\}$ is detected in $\K(\K(R))$. 
\end{cor}
\begin{proof}
Since there exists a map of $E_2$ ring spectra $R \rightarrow H\mathbb{F}_q$ inducing a map $\K(\K(R))\rightarrow \K(\K(\mathbb{F}_q))$. Therefore, there is a commutative diagram 
\begin{align}\label{diag ring of int}
	\xymatrix{ 
		S \ar[r] & \K(\K(S))  \ar[r]  \ar[dr] &  \K(\K(R)) \ar[d] \\ 
						& & \K(\K(\mathbb{F}_q)).
						}
\end{align}
Since the $\beta$-family is non-trivial in the image of the map $\pi_*S\to \pi_*\K(\K(\mathbb{F}_q))$, it is also nontrivial in the image of the unit map $\pi_*S\to \pi_*(\K(\K(R)))$. 
\end{proof}
In particular, this implies the following result:
\begin{cor}\label{last cor}
The $p$-primary beta family $\{ \beta_s : s \not\equiv 0 \mod p\}$ is detected in $\K(\K(\mathbb{Z}))$ for all $p\ge 5$.  
\end{cor}
\begin{rem}
Note that the $\alpha$-family is detected in $\K(\bZ)$. Since $\K_0(\mathbb{Z})\cong \mathbb{Z}$, there is a map of commutative ring spectra $\K(\mathbb{Z})\to H\mathbb{Z}$. We may consider the infinite family of maps 
\[ S\to \ldots \to \K(\K(\K(\bZ)))\to \K(\K(\bZ))\to \K(\bZ) \]
and a specialization of the Greek-letter family red-shift conjecture is that the $p$-primary $n$-th Greek letter family is in the image of the unit map $S\to \K^{(n)}(\bZ)$, where $\K^{(n)}(\bZ)$ is algebraic K-theory iterated $n$-times, for $p$ sufficiently large. Corollary \ref{last cor} is therefore evidence for this conjecture when $n=2$. 
\end{rem} 
\bibliography{sources}
\bibliographystyle{alpha}
\end{document}